\documentclass[a4paper, 11pt]{article}

\pdfoutput=1

\usepackage{graphicx}

\usepackage{color}
\usepackage{amsmath}
\usepackage{amsfonts}
\usepackage{amssymb}
\usepackage{amsbsy}
\usepackage{amsthm}
\usepackage{epstopdf}
\usepackage{algorithm}
\usepackage{setspace}
\usepackage[cm]{fullpage}
\usepackage{url}
\usepackage{multirow}

\theoremstyle{remark}
\newtheorem{remark}{Remark}
\newtheorem{problem}{Problem}

\newtheorem{theorem}{Theorem}

\newcommand{\T}{\ensuremath{\mathcal{T}} }
\newcommand{\E}{\ensuremath{\mathcal{E}} }
\newcommand{\NN}{\ensuremath{\mathcal{N}} }

\newcommand{\intO}{\int_\Omega}
\newcommand{\R}{\mathbb{R}}
\def\Rd{{\mathbb{R}^d}}
\def\vz{v^*}
\def\dvg{{\rm div}}
\def\ben{\begin{eqnarray}}
\def\een{\end{eqnarray}}
\newcommand{\ys}{q^*}

\newcommand{\Lam}{\Lambda}

\DeclareMathOperator*{\argmin}{argmin}

\renewcommand{\L}{l}
\newcommand{\Lspace}{\ensuremath{\mathcal{L}} }

\title{A FEM approximation of  a two-phase obstacle problem \\
and its %functional
a posteriori error estimate}

\author{FARID BOZORGNIA$^{*}$ and JAN VALDMAN$^{**}$\\ \\
\noindent
$^*$ Department of Mathematics, Instituto Superior T\'{e}cnico, \\
Av. Rovisco Pais, P-1049-001 Lisboa, Portugal. \\ \\
$^{**}$ Institute of Mathematics and Biomathematics, Faculty of Science, \\
University of South Bohemia, Brani\v sovsk\' a 31, CZ--37005, Czech Republic \\
and Institute of Information Theory and Automation, Academy of Sciences, \\
Pod vod\'{a}renskou v\v{e}\v{z}\'{\i}~4, CZ--18208~Praha~8, Czech Republic.
}

\begin{document}

\maketitle

\begin{abstract}
This  paper is concerned with the two--phase obstacle problem, a type of a variational free boundary problem. We recall the basic  estimates of \cite{ReVa2015} and verify them numerically on two examples in two space dimensions. A solution algorithm is proposed for the construction of the finite element approximation to the two--phase obstacle problem. The algorithm is not based on the primal (convex and nondifferentiable)  energy minimization problem but on a dual maximization problem formulated for Lagrange multipliers. The dual problem is  equivalent to a quadratic programming problem with box constraints.
The quality of approximations is measured by a functional a posteriori error estimate which provides a guaranteed upper bound of the difference of approximated and exact energies of the primal minimization problem. The majorant functional in the upper bound contains auxiliary variables and it is optimized with respect to them to provide a sharp upper bound.
A space density of the nonlinear related part of the majorant functional serves as an indicator of the free boundary. %The value of the optimized majorant functional gives automatically a lower bound of the exact energy.

%If a reference solution is available, a lower bound of the difference of approximated and exact energies of the primal minimization problem is also evaluated.
\end{abstract}

\section{Introduction}

A free boundary problem is a partial differential equation where the equation changes qualitatively across a level set of the equation solution $u$ so the part of the  domain where   the equation changes is a priori  unknown.
A general form of  elliptic  free boundary problems    can be written as
\begin{equation}\label{G}
\Delta u= f(x, u,\nabla u) \qquad \text{ in } \Omega,
\end{equation}
where  the right hand side term is  piecewise continuous, having jumps at some values of the arguments $u$ and $\nabla u.$ Here  $\Omega$ is a bounded open subset of $\mathbb{R}^n$ with smooth  boundary   and   Dirichlet boundary conditions are considered. In this paper we are concerned about the particular elliptic  free boundary problem
\begin{equation}\label{H}
\left \{
\begin{array}{ll}
\Delta u=  \alpha_{+}  \chi_{\{u >0\}}-\alpha_{-} \chi_{\{u <0\}} &  \text{ in } \Omega,\\
   u=g     &   \text{ on }  \partial \Omega.
 \end{array}
\right.
  \end{equation}
Here, $\chi_{A}$ denotes the  characteristic function of the set  $A$,  $\alpha_{\pm}:\Omega\rightarrow \mathbb{R}$ are positive  and Lipschitz continuous functions and  $ g \in W^{1,2}(\Omega) \cap L^{\infty}(\Omega)$ and $g$ changes  sign on  $ \partial \Omega$. The boundary
\[
 \left( \partial \{x\in \Omega: u(x)>0 \}  \cup \partial\{x\in \Omega: u(x) < 0\} \right) \cap \Omega,
 \]
 is called the free boundary. Properties  of the solution of  the two-phase obstacle problem, regularity   of solution and  free boundary have been studied in \cite{Ur, We}. Its is known that the differential equations from \eqref{H} represents the Euler-Lagrange equation corresponding  to the minimizer of  the functional
\begin{equation}\label{eq:energy}
J(v)=\int_{\Omega}\left(\frac{1}{2}|\nabla v|^2+\alpha_{+} \text {max}(v,0) + \alpha_{-} \text{max}(-v,0)\right)dx
\end{equation}
over the affine space  
\begin{equation} 
K=\{ v \in W^{1,2}(\Omega): v-g \in W^{1,2}_{0}(\Omega) \}.
\end{equation}
The functional $J: K \rightarrow \mathbb{R}$  is convex, coercive on $K$ and weakly lower semi-continuous, hence  the minimum of $J$ is attained  at some  $u \in  K$. The following minimization problem is therefore uniquely solvable.

 \begin{problem}[Primal problem]
\label{prob1} Find $u \in K$ such that
\begin{equation}
J(u)=\inf\limits_{v \in K} J(v).   \label{eq:energy_inf}
\end{equation}
\end{problem}
Note that if we let $ \alpha^{-} =0, $ and assume  that $g$ is  nonnegative  on the boundary, then  we obtain the well-known  one-phase obstacle problem,  see e.g. \cite{Caf,KS}.

There are  numerous  papers on approximations and error analysis for  the one-phase obstacle problem in terms on variational inequalities \cite{FALK, GLT}. In \cite{No}  a sharp $L^{\infty}$ error estimate for semilinear elliptic problems with free boundaries is given. For  obstacle problem and combustion problems, the author  uses  regularization  of penalty term combined with piecewise linear finite elements  on a triangulation and then shows that the method is accurate  in $L^{\infty}$. Using  non-degeneracy property of one-phase obstacle problem, a sharp interface error estimate is derived.   In \cite{NOS}  error estimates for  the finite element approximation of the  solution  and free boundary of the   obstacle problem   are presented.   Also an optimal error analysis for the thin obstacle problem is derived.

Recently,  the numerical approximation of  the two-phase membrane problem has  attracted much interests. Most approximations are based on  the finite difference methods. In \cite{Bo} different methods  to approximate the solution are  presented. The first method  is based on properties of  the given free boundary problem and exploit the disjointness of positive and negative parts of the solution. Regularization method and   error estimates  are  given. The  a priori   error   gives a computable estimate  for gradient of the error  for  regularized  solutions in   the  $L^{2}$. In \cite{ARM}, the authors   rewrite the two phase obstacle problem in an equivalent  min-max  formula then for this new form they  introduce
 the notion of viscosity solution.  Discritization of the min-max formula  yields   a certain linear approximation system. The existence and uniqueness of the solution of the discrete nonlinear system  are shown.  Also in \cite{avetikparab} the author presents a finite difference approximation for a parabolic version of  the two-phase membrane problem.

A finite element  scheme for solving obstacle problems in divergence form is given in \cite{TSFO}.  The authors  reformulate the obstacle in terms of an  $L^{1}$ penalty on the variational problem.  The reformulation is an exact regularizer in the sense that for large  penalty parameter, it can recover the exact solution.   They  applied the scheme   to approximate classical elliptic obstacle problems, the two-phase membrane problem and the Hele-Shaw model.

We propose a different finite element scheme for solving the two-phase obstacle problem based on the dual maximization problem for Lagrange multipliers. The main focus of the paper is  the verification of a posteriori error estimates developed in \cite{ReVa2015}. For any obtained FEM approximation $v$, we can explicitely compute the upper bound of the difference $J(v)-J(u)$ of the approximate  energy $J(v)$ and of the exact unknown minimal energy $J(u)$. Since this upper bound is quaranteed we automatically have a lower bound of  the exact energy $J(u)$.  The studied aposteriori error estimates also provide the approximate indication of the exact free boundary. This is demontrated on two numerical tests in two space dimentions. A MATLAB code is freely available for own testing.  

The structure of paper  is as follows. In Section 2,  we  present an overview of  basic concepts and mathematical background and recall energy and majorant estimates of \cite{ReVa2015}. Section 3 deals with discretization using finite elements: construction of the FEM approximation (Algorithm 1) and the optimization of the functional majorant (Algorithm 2). Section 4 reports on numerical examples and section 5 concludes the work.

\section{Mathematical Background and Estimates}
Elements of convex analysis are used throughout this paper, in particularly the duality method by conjugate functions  \cite{ET} .
For reader's convenience, let us summarize
the basic notation used in what follows:

\bigskip

%\vspace{.7em}

\hspace*{-1.6em}\fbox{
\begin{minipage}[t]{0.46\linewidth}
\small

$n=1, 2, 3$ dimension of the problem,

$\alpha_{\pm} \geq 0$ problem coefficients,

$v^{\pm}$  positive and negative parts of the function,

$J: K \rightarrow \R$ primal functional to be minimized,

$J^{*}: Q^* \rightarrow \R$ conjugate functional to be maximized,

$J_{\mu} (\cdot)$ perturbed functional with multiplier $\mu$,

$u,v \in K$ exact and arbitrary minimizers of  $J$,

$p^{*}, q^{*} \in Q^{*}$ exact and arbitrary maximizers of  $J^{*}$,

$\lambda, \mu \in \Lambda$ exact and arbitrary multipliers,

$D_J(\cdot, \cdot): K \times Q^* \rightarrow \R$ compound functional,

$M_+: K \times \R_+ \times Y^* \times \Lambda  \rightarrow \R$ majorant functional,

$\eta^* \in Y^*$ flux variable in $M_+$ approximating $p^{*}$,

$C_\Omega$ constant from generalized Friedrich's inquality,

\end{minipage}\ \

\begin{minipage}[t]{0.48\linewidth}\small

$\mathcal T_h$ uniform regular triangular mesh with mesh size $h$,
 $\NN_I, \NN_D$ internal and Dirichlet nodes, $I, D$ their indices,
$|\E|, \, |\NN|, \, |\T|$ number of edges, of nodes and of triangles,

$K_h, Q^*_h, \Lambda^*_h$ finite element approximation spaces on $\mathcal T_h$

$I^*:  \Lambda^*_h  \rightarrow \R$ discrete dual energy to be maximized,

$\mathbb{K}$ stiffness matrix in $K_h$,

$\mathbb{K}_{I,I}, \mathbb{K}_{I,D}, \mathbb{K}_{D,D}$ its subblocks with respect to $I$ and $D$,

$\mathbb{M}$ generalized mass matrix ($L^2$ - product of $K_h$ and $\Lambda_h$),

$\mathbb{M}_I, \mathbb{M}_D$ its subblocks with respect to to $I$ and $D$,

 $\boldsymbol{\lambda}, \boldsymbol{\mu}, \textbf{v},  \textbf{u}_{\boldsymbol{\lambda}}$ discrete vectors,

  $ \textbf{v}_I, \textbf{v}_D, \textbf{u}_{{\boldsymbol{\lambda}}_I}, \textbf{u}_{{\boldsymbol{\lambda}}_D}$ subvectors with respect to $I$ and $D$,

$u_{ref}$ reference solution,

$M_{+1}, M_{+2}, M_{+3}$ majorant functional subparts,

\end{minipage}\medskip
}

\vspace{.3em}

\begin{center}
{\small\sl List \,1.\ }
%\begin{minipage}[t]{.7\textwidth}\baselineskip=8pt
{\small
Summary of the basic notation used thorough out this  paper.
}
%\end{minipage}
\end{center}

% We refer to
%\cite{Han} for the general relations and to \cite{NeRe, ReGruyter} for detailed forms for conjugate and compound fucntionals in the next section. Two main %presented estimates, energy and majorant estimate were already derived in \cite{ReVa2015} and the meaning of this section is to show them in more %general context. \\
Let $V$ and $Q$  be two normed spaces, $ V^{*}$ and $Q^{*}$ their dual spaces and let  $\langle \cdot,\cdot \rangle $ denote the duality pairing.
 Assume that there exists   a continuous  linear operator $\L$ from $V $ to $Q$,
$  \L \in {\Lspace}(V,Q). $
 The adjoint operator $\L^{*} \in {\Lspace}(Q^{*}, V^{*}) $ of the operator $\L$ is defined through  the relation
\[
\langle \L^{*}q^{*} ,v \rangle=\langle q^{*}, \L v \rangle \ \ \  \forall v \in V, \ q^{*} \in  Q^{*}.
\]
Let $J: V \times Q \rightarrow \mathbb{\overline{R}}$ is a convex functional mapping in the space of extended reals $\mathbb{\overline{R}}=\mathbb{R} \cup {\{ -\infty,+\infty }\}$. 
Consider the minimization problem
\begin{equation}\label{M}
\underset{v \in V}   \inf  J(v, \L v).
\end{equation}
and its dual conjungate problem
\begin{equation}\label{MM}
 \underset{q^{*} \in Q^{*}} \sup[ -J^{*}(\L^{*}q^{*},-q^{*})],
 \end{equation}
where the  convex conjugate function of $J$ is given by
\[
J^{*}(v^{*},q^{*})=\underset{v \in V,\,q \in Q} \sup [\langle v,v^{*}\rangle + \langle q, q^{*}\rangle - J(v,q)], \ v^{*} \in  V,\, q^{*} \in Q^{*}.
\]
The relation between (\ref{M}) and (\ref{MM})  is stated in the following theorem.
\begin{theorem}[Theorem 2.38 of \cite{Han}]
\label{relation}
Assume that $V $ is a reflexive Banach space and $Q$  is a normed  vector space, and let  $\L \in {\Lspace}(V,Q). $
 Let  $J: V \times Q \mapsto \mathbb{\overline{R}}  $ be  proper lower semi continuous,   strictly convex such that
 \begin{enumerate}
    \item There exists  $v_{0}\in V,$ such that $  J(v_{0}, \L v_{0})< \infty$ and $ q \rightarrow J (v_{0},q) $ is continuous at $ \L v_{0}.$
  \item $J(v, \L v) \rightarrow +\infty,$ as $ \|v\|\rightarrow \infty, v \in V.$
   \end{enumerate}
    Then  problem (\ref{M}) has a solution $ u \in V $ also  problem (\ref{MM}) has a solution   $ p^{*} \in Q^{*}$, and
 \begin{equation}\label{ajab}
 J(u,\L u)= -J^{*}(\L^{*}p^{*},-p^{*}).
\end{equation}
\end{theorem}
In the case that  the function $J$ is of a separated form, i.e.,
\[
J(v,q)=F(v)+G(q)  \qquad v \in V , q \in Q,
\]
then the conjugate of $J $ is
\[
J^{*}(v^{*},q^{*})=F^{*}(v^{*})+G^{*}(q^{*}),
\]
where $F^{*}$ and $G^{*}$      are the conjugate functions of $F $ and $G$, respectively.
 To calculate the conjugate function when  the functional  is defined by an integral,  we use the following theorems which can be found in \cite{Han}.
\begin{theorem}[Theorem 2.35 of \cite{Han}]
Assume $ h: \Omega \times \mathbb{R}^n \longrightarrow \mathbb{R}$ is a Carath\'{e}odory function with  $ h \in L^{1}(\Omega)$  and suppose
\[
G(q)=\int_{\Omega} h(x,q(x))\, dx.
\]
Then the conjugate function of $G$ is
\[
G^{*}(q^{*})=\int_{\Omega} h^{*}(x,q^{*}(x)) dx \  \ \ \forall q^{*} \in Q^{*},
\]
where
\[
h^{*}(x,y)= \underset{\xi \in  \mathbb{R}^n }  \sup[y \cdot \xi-h(x,\xi)].
\]
\end{theorem}
%
%
%Let $ u \in V $ be a solution of  the minimization  problem (\ref{M}) and $q^{*}$ be defined as in (\ref{ajab}). For any $ v\in V, $  define the energy difference
%\[
%ED(u,v)= J(v,\textrm{A}v)- J(u,\textrm{A}u).
%\]
%\begin{theorem}\label{E}
%Suppose the assumptions in  Theorem \ref{relation}  are satisfied. Then
%\[
%ED(u,v)\leq  J(v,\textrm{A}v)+ J^{*}(A^{*}q^{*},-q^{*}) \ \ \ \forall v \in V, q^{*} \in Q^{*}.
%\]
%\end{theorem}
%

The compound functional $D_{J}(v,q; v^*, q^{*}): (V \times Q) \times (V^* \times Q^*) \rightarrow \mathbb{\overline{R}}$  is defined by
\begin{equation}
D_{J}(v,q; v^*, q^{*}):= J(v, q)+ J^*(v^*, q^{*}) - \langle v, v^{*}  \rangle - \langle q, q^{*}  \rangle .
\end{equation}
It holds $D_{J}(v,q; v^*, q^{*}) \geq 0$ for all $(v, q) \in V \times Q$, $(v^*, q^*) \in V^* \times Q^*$ and  $D_{J}(v,q; v^*, q^{*})=0$ only if the function $(v, q)$ belongs to set of subdifferential  $ \partial J^{*}(v^*, q^{*})$ and $(v^*, q^{*}) $ belongs to  set of subdifferential  $ \partial J(v,q),$  see Proposition 1.2 of \cite{ReGruyter}. 

\subsection{Energy identity}
%We are ready to derive a  posteriori error estimate of the functional error type   for Problem \ref{prob1}.
%  To see more about derivation we refer to  \cite{ReVa2015} and its derivation is based on general results for variational inequalities from Section 7 in \cite{NeRe}.

For simplicity of notation, we introduce the  positive and negative parts of a function $v$
\begin{equation*}
 v^{+} :=  \text {max}(v,0), \qquad    v^{-} :=  \text {max}(-v,0),    
\end{equation*}
so it holds  $v=v^{+}-v^{-}$ and $|v|=v^{+} + v^{-}$. The Euclidean norm in   $\mathbb{R}^{n}$ is denoted by  $|\cdot|.$  We write \eqref{eq:energy}  in the form
$ J(v) = F(v) + G(\L v),$
where
\begin{equation}
\begin{aligned}
\label{GF}
F(v):&=
   \intO \Big( \alpha_{+}  v^{+} + \alpha_{-}  v^{-} \Big) \,   dx, \qquad 
 G(\L v):&=\frac12 \intO   \nabla v\cdot\nabla v  \,dx \\  
\end{aligned}
\end{equation}
and $\L: K \rightarrow Q=L^2(\Omega,\Rd)$ is the gradient operator $\L v = \nabla v$.

\begin{remark}
Theorem 1 assumes $J:  V \times Q \rightarrow \mathbb{\overline{R}} $, where $V$ is a normed space. It can be shown that all results are also valid for $J:  K \times Q \rightarrow \mathbb{R} $ from above. 
\end{remark}

The corresponding dual conjugate functionals are
\begin{equation}\label{GF_duals}
\begin{aligned}
F^*(v^*)&= \intO  h^* v^*   \,dx, \qquad 
G^*(\ys)&= \frac12 \intO  \ys \cdot \ys  \,dx, 
\end{aligned}
\end{equation}
where $h^*(z^*)=0$ for $z^* \in [- \alpha_{-}, \alpha_{+} ]$ otherwise $h^*(z^*)=+\infty$.
%Using $G^*(\ys)$ we define the Lagrangian
%$$L(v,\ys):=  %\intO  \ys \cdot \Lambda v  \,dx
%\left<\ys, \L v   \right> -  G^*(\ys) + F(v), \quad   \ys \in Y^{*}$$
%and show relations
%\begin{equation}
%\begin{aligned} \label{duality}
%J(u)&=\inf\limits_{v \in K} J(v) = \inf\limits_{v \in K} \sup\limits_{\ys \in \ys} L(v,\ys) = \sup\limits_{\ys \in \ys} \inf\limits_{v \in K}  L(v,\ys)  \\
%&= \sup\limits_{\ys \in \ys} \left( -  G^*(\ys) + \inf\limits_{v \in K}  \left(  \left<\ys, \L v   \right> +   F(v)  \right) \right) \\
%&= \sup\limits_{\ys \in \ys} \left( -  G^*(\ys) - \sup\limits_{v \in K}  \left(  \left<-\L^* \ys, v   \right> -   F(v)  \right) \right) \\
%&= \sup\limits_{\ys \in \ys} \left( -  G^*(\ys) - F^*(-\L^* \ys )  \right)  =: \sup\limits_{\ys \in \ys}  -J^{*}(\L^{*}q^{*},-q^{*}) \\
%&=: \sup\limits_{\ys \in \ys} I^*(\ys) =  I^*(p^*) ,
%\end{aligned}
%\end{equation}
%where $I^*(\cdot)$ denotes a dual functional   {\color{red} Is  $I^*(\cdot)$ the dual of  $J(u)$?}
%Appying \eqref{relation}, we immediately obtain
%$$G(\L u) + F(u) =\inf\limits_{v \in K} \left( G(\L v) + F(v) \right)= \sup\limits_{\ys \in \ys} \left( -  G^*(\ys) - F^*(-\L^* \ys )  \right) = \sup\limits_{\ys \in \ys}  -J^{*}(\L^{*}q^{*},-q^{*}) .$$
%It can be show that the supremal value on the right side is obtained for
Since $\L v =\nabla v$, the dual operator  is represented by the divergence operator $-\L^* \ys = \dvg  \ys $.
Combining \eqref{GF} and \eqref{GF_duals} we derive compound functionals
\begin{eqnarray}
D_F(v,\vz)=\intO\left(\alpha_{+} v^{+} + \alpha_{-}  v^{-}- \vz\,  v \right) \,dx, \qquad 
D_G(\L v,\ys)
=\frac12\intO (\L v -  \ys)  \cdot (\L v  -  \ys)\, dx,  \label{DG_compound} 
\end{eqnarray}
where the form for $D_F(\cdot)$ is valid if the condition
\ben
v^*  \in [- \alpha_{-}, \alpha_{+} ]
\een
is satisfied almost everywhere in $\Omega$, otherwise $D_F(v,\vz)=+\infty$. For the gradient type problem, it holds
\ben \label{flux_exact_choice}
p^*=\nabla u,
\een
i.e., $p^*$ represents the exact flux (gradient of the exact solution). Compound functionals appear in the {\it energy identity}  (Proposition (7.2.13) of \cite{NeRe})
\ben
\label{mainidentity}
D_F(v,-\L^*p^*) + D_G(\L v,p^*)=J(v)-J(u)  \qquad \mbox{ for all } v \in K.
\een

The terms in the left part of \eqref{mainidentity} are
\begin{eqnarray}
\label{compound_F}
 D_F(v,-\L^*p^*)&=&\intO \left(\alpha_{+}   v^{+} + \alpha_{-}  v^{-}-( \dvg p^*)v
\right)dx, \\
\label{compound_G_v}
D_G(\L v,p^*)
%=\intO\left(\frac12 A\nabla v\cdot\nablav+\frac12 A^{-1}p^*\cdot p^*-\nabla v\cdot p^*\right)dx
&=&
\frac12\intO \nabla (u-v)\cdot\nabla (u-v)dx=\frac12\| \nabla(u-v) \|^{2}_{L^{2}(\Omega)}
\end{eqnarray}
and $D_F(v,\vz)$ is always finite since the condition $\dvg p^{*}   \in [ -\alpha_{-}, \alpha_{+} ]$ is always satisfied.

\begin{remark}[Gap in the energy estimate]
If we drop the nonnegative term
$D_F(v,\vz)$ we get
\ben
\label{inequality_energy}
\frac12\| \nabla(u-v) \|^{2}_{L^{2}(\Omega)} \leq J(v)-J(u) \qquad \mbox{ for all } v \in K
\een
which is well known in connection to class of nonlinear problems related to variational inequalities. For the two-phase obstacle problem, is was derived in \cite{Bo}.  The gap in the sharpness of the estimate \eqref{inequality_energy}  is exactly measured by the term $D_F(v,-\L^*p^*)$. By respecting $D_F(v,-\L^*p^*)$ we can get the equality formulated by the main estimate \eqref{mainidentity}. The contribution of $D_F(v,-\L^*p^*)$ is expected not to be very high for good quality approximation $v \in K$ to the exact solution $u$. An example, when the gap becomes significantly large (for a bad approximation $v$) is given in Section 4 of \cite{ReVa2015}.  The form of $D_F(v,-\L^*p^*)$ represents a  certain measure of the error associated with free boundary and it is further discussed in Section 2 of  \cite{ReVa2015}.
\end{remark}

\subsection{Majorant estimate}\label{2.2}

The exact energy $J(u)$ in the energy identity \eqref{mainidentity} and the energy inequality \eqref{inequality_energy} is not computable without the knowledge of the exact solution $u$. However, we can get its computable lower bound using a perturbed functional
\begin{equation}\label{def:perturbed_functional}
J_{\mu}(v) : = G(\L v) + F_{\mu}(v), \qquad F_{\mu}(v) :=   \intO   \mu  v  \,dx,
\end{equation}
where a multiplier $\mu \in \Lam$ belongs to the space
$$ \Lam:=\left\{\mu \in L^\infty(\Omega):
\mu(x) \in [-\alpha_{-}, \alpha_{+} ] \; \mbox{a.e.}\; \textrm{in}\; \Omega\right\}.$$
The perturbed functional $J_{\mu}(v)$ replaces the non-differentiable functional $J(v)$ at the cost of a new variable $\mu \in \Lam $ in $F_{\mu}(\cdot)$.
It holds
\begin{equation}\label{inequalities}
J(u) = \inf\limits_{v\in K} \sup\limits_{\mu \in \Lam} J_{\mu}(v) =
\sup\limits_{\mu \in \Lam} \inf\limits_{v\in K} J_{\mu}(v) \geq
\inf\limits_{v\in K} J_{\mu}(v)  =: J_\mu(u_{\mu} )    \quad \mbox{for all } \mu \in \Lam,
\end{equation}
 where $u_{\mu} \in K$ is unique. In view of \eqref{inequalities}, the minimal perturbed energy $J_\mu(u_{\mu} )$ serves as the lower bound of  $J(u)$.
%In applications, we are usually able to compute a discretization of $J_\mu(u_{\mu} )$ and therefore only its upper bound.
We find a computable lower bound of $J_\mu(u_{\mu} )$ by means of the dual counterpart of the perturbed problem. The dual problem is generated
by  the Lagrangian
$$L_{\mu}(v,\ys):=  %\intO  \ys \cdot \L v  \,dx
\left<\ys, \L v   \right> -  G^*(\ys) + F_{\mu} (v).$$
We note $v=g+w,$ where $w \in H^1_0(\Omega)$ and estimate
\begin{eqnarray}
J_\mu(u_{\mu} )&=&\inf_{v \in K} J_\mu(v)= \inf_{v \in K}  \sup_{q^* \in Q^*}  L_\mu(v,\ys) = \sup_{\ys \in Q^* }  \inf_{v \in K}    L_\mu(v,\ys) \nonumber \\
 &=& \sup_{\ys \in Q^* }   \inf_{w \in H^1_0} \{  -  G^*(\ys)  +   \intO (\ys \cdot  \nabla g -  \mu g) \,\textrm{d}x  +  \intO (\ys \cdot \nabla w - \mu w ) \,\textrm{d}x \} \nonumber \\
&\geq & \sup_{\ys \in  Q_{\mu} }  \{  -  G^*(\ys)  +   \intO (\ys \cdot  \nabla g -  \mu g) \,\textrm{d}x   \} \nonumber \\
&=&  \sup_{\ys \in Q_{\mu}} J_\mu^*(\ys) \geq  J_\mu^*(\ys) \qquad \mbox{for all }  \ys \in Q^{*}_{\mu}, \label{eq:inf_sup}
\end{eqnarray}
where
\begin{equation}\label{eq:energy_dual}
J_\mu^*(\ys): =  -  G^*(\ys) + \intO (\ys  \cdot \nabla g -  \mu g) \,\textrm{d}x
\end{equation}
and
\begin{equation}\label{def:Qfmu}
Q^*_{\mu} = \{ \ys \in  Q^*:  \intO (\ys \cdot \nabla w - \mu w ) \,\textrm{d}x =0 \mbox{ for all } w \in H^1_0(\Omega)  \}.
\end{equation}

Due to  \eqref{inequalities} and  \eqref{eq:inf_sup}, we obtain the estimate
\begin{equation}\label{ineq:energy_perturbed}
 J(v) - J(u)  \leq J(v) - J_\mu^*(\ys) = \frac{1}{2} ||  \nabla v - \ys ||_{L^{2}(\Omega)}^2 +
\intO \Big( \alpha_{+}  v^{+} + \alpha_{-}   v^{-}  -  \mu  v \Big)\,\textrm{d}x
\end{equation}
 valid for all $v \in K, \mu \in \Lam, \ys \in Q^*_{\mu}$. The right-hand size of \eqref{ineq:energy_perturbed} is fully computable, but it requires the constraint $\ys \in Q^*_{\mu}$. To bypass this constraint we  introduce  a new variable
 $$\eta^* \in Y^*: = H(\Omega, \rm{div})$$
 and project it to $Q^*_{\mu}$. The space $ H(\Omega, \rm{div})$ is a subspace of $ L^2(\Omega, \R^n)$ that contains vector-valued functions with square-summable divergence.
There holds a projection-type inequality
(see, e.g.,  Chapter 3 of \cite{ReGruyter})
$$
\inf_{\ys \in Q_{\mu}}  ||\eta^* - \ys ||_{L^{2}(\Omega)}  \leq C_{\Omega} || {\rm{div}} \,\eta^{*} + \mu||_{L^{2}(\Omega)} \qquad \mbox{for all } \eta^* \in Y^*,
$$
where the constant $C_{\Omega}>0$ originates from the generalized Friedrichs inequality
$$|| w ||_{L^{2}(\Omega)} \leq C_{\Omega} ||\nabla w ||_{L^{2}(\Omega)} \qquad \mbox{for all } w \in H^{1}_{0}(\Omega). $$ Then, the $\ys$-dependent term in \eqref{ineq:energy_perturbed}  satisfies
\begin{eqnarray}
\inf_{\ys \in Q^*_{\mu}} ||  \nabla v - \ys ||_{\Omega}^2 &\leq &
  ( ||   \nabla v - \eta^* ||_{\Omega} + \inf_{\ys \in Q^*_{\mu}}  ||\eta^* - \ys ||_{\Omega} )^2 \nonumber \\
&\leq&(     ||  \nabla v - \eta^* ||_{\Omega} + C_{\Omega} || {\rm{div}} \,\eta^{*} + \mu||_{\Omega} )^2 \label{ineq:unconstrained} \\
&\leq&  (1+\beta) ||  \nabla v - \eta^* ||_{\Omega}^2+ (1+\frac{1}{\beta})  C_{\Omega}^2 || {\rm{div}} \,\eta^{*} + \mu||_{\Omega}^2, \nonumber
\end{eqnarray}
where we used  Young's estimate with a parameter $\beta > 0$ in the last inequality. Hence the combination of \eqref{ineq:energy_perturbed} and  \eqref{ineq:unconstrained} yields {\em the majorant estimate}
\begin{eqnarray}\label{eq:majorant_estimate}
J(v) - J(u) \leq  M_+(v; \beta, \eta^{*}, \mu  ) \qquad \mbox{ for all } v \in K, \mu \in \Lam, \eta^* \in Y^*, \beta>0,
\end{eqnarray}
where a nonnegative functional
\begin{multline}  \label{majorant}
M_+(v; \beta, \eta^{*}, \mu   ) :=  \frac{1}{2}     (1+\beta) ||  \nabla v - \eta^* ||_{\Omega}^2 \\
+ \frac{1}{2}   (1+\frac{1}{\beta})  C_{\Omega}^2 || {\rm{div}} \,\eta^{*} - \mu||_{\Omega}^2  \\  + \intO \Big( \alpha_{+}  v^{+} + \alpha_{-}   v^{-}  -  \mu  v \Big)\,\textrm{d}x,
\end{multline}
represents a functional error majorant.

\begin{remark}
The final form of the functional error majorant  \eqref{majorant} is slightly different to formula (3.13) in  \cite{ReVa2015}. In \eqref{majorant}, only one multiplier variable $\mu$ is introduced replacing  two multipliers $\mu_{-}, \mu_{+}$ of \cite{ReVa2015}. Another simplification in this paper is that the variable diffusions coefficient matrix $A$ is not considered here and we treat the quadratic part of the energy $\frac12 \intO   \nabla v\cdot\nabla v  \,dx$  instead of $\frac12 \intO   A \nabla v\cdot\nabla v  \,dx$.

\end{remark}

\section{Discretization}

We assume a domain $\Omega \subset\mathbb{R}^2$ with a polygonal boundary  discretized by an uniform regular triangular mesh $\mathcal T_h$  in the sense of Ciarlet \cite{CiaBook}, where $h$ denotes the mesh size.  Let $\E$ denote the set of all edges and $\NN$ the set of all nodes in
$\T_h$.  By
$$|\E|, \, |\NN|, \, |\T|$$ we mean  the number of edges, of nodes and of triangles of $\T_h$.
%A generalization to the case three-dimentional case  $\Omega \subset\mathbb{R}^3$  is then straightforward.
The following lowest order finite elements (FE) approximations are considered:
\begin{itemize}
\item The exact solution $u \in K$ of the  two phase  obstacle problem is approximated by
$$u_h \in K_h: = K \cap P_1(\mathcal T_h),$$ where  $P_1(\mathcal T_h)$ denotes the space of elementwise nodal and continuous functions defined on $\mathcal T_h$.
\item The exact multiplier $\lambda \in \Lam$ is approximated by
$$\lambda_h \in \Lam_h := \Lam \cap P_0(\mathcal T_h), $$ where $P_0(\mathcal T_h)$ denotes the space of element wise constant functions defined on $\mathcal T_h$.
\item The  flux variable $\eta^* \in Q$ in the functional majorant is approximated by
$$\eta^{*}_h \in Y^*_h:= RT_0(\mathcal T_h), $$
where $RT_0(\mathcal T_h)$ is the space of the lowest order Raviart-Thomas functions.
\end{itemize}
Note that dimensions of these approximation spaces are 
$$\dim(K_h)=|\NN|, \quad \dim( \Lam_h)=|\T|, \quad \dim( Y^*_h)=|\E|.$$

We are  interested in two computation tasks:  First obtaining a discrete solution $u_h$ or its approximation $v_h$ and then  measuring its quality by  the optimized functional error majorant.

\subsection{Dual based algorithm for Lagrange multipliers}\label{subsec:solver}
 As we mentioned in \ref{2.2}  due to the  non-differentiability  term in  $J(\cdot)$, we do not solve the approximative solution $u_h \in K_h$  from the relation $J(u_h) = \inf\limits_{v\in K_h}  J(v)$ directly. We estimate
\begin{equation}\label{inequalities_discrete}
J(u_h) = \inf\limits_{v_h \in K_h} \sup\limits_{\mu \in \Lam} J_{\mu}(v) \geq  \inf\limits_{v_h \in K_h} \sup\limits_{\mu_h \in \Lam_h} J_{\mu_h}(v_h) \geq
\sup\limits_{\mu_h \in \Lam_h} \inf\limits_{v_h \in K_h} J_{\mu_h}(v_h)  %= \sup\limits_{\mu_h \in \Lambda_h} I(\mu_h) =
=:J_{\lambda_h}(u_{\lambda_h})
\end{equation}
%where $ I(\mu_h):=  \inf\limits_{v_h \in K_h} J_{\mu_h}(v_h).$
and look for an approximation pair $(\lambda_h, u_{\lambda_h}) \in \Lam_h \times K_h$ instead.
Note, in general $ u_{\lambda_h} \not= u_h$ and it holds $J(u_{\lambda_h}) \geq J(u_h)$ only.
The  saddle point problem on the right-hand side of \eqref{inequalities_discrete} can be  further reformulated as a dual problem for a Lagrange multiplier
\begin{equation}\label{dual_discrete}
I^*(\lambda_h) = \sup\limits_{\mu_h \in  \Lam_h} I^*(\mu_h), \qquad \mbox{where } I^*(\mu_h) := \inf\limits_{v_h \in K_h} J_{\mu_h}(v_h).
\end{equation}

Approximations  $v_h \in K_h$ and $\mu_h \in \Lam_h$ from \eqref{def:perturbed_functional}  are equivalently represented by discrete column vectors $\textbf{v} \in \mathbb{R}^{|\NN|}, \boldsymbol{\mu} \in \mathbb{R}^{|\E|} $ and we can rewrite $J_{\mu_h}(\cdot)$  from \eqref{def:perturbed_functional}  as
\begin{equation}\label{discrete_J_mu}
J_{\boldsymbol{\mu}}(\textbf{v}) = \frac{1}{2} \textbf{v}^T \mathbb{K} \textbf{v} +  \textbf{v}^T \mathbb{M} \boldsymbol{\mu},
\end{equation}
where $ \mathbb{K}  \in \mathbb{R}^{|\NN| \times |\NN|} $ and $\mathbb{M} \in \mathbb{R}^{|\NN| \times |\T|} $. The square matrix  $ \mathbb{K}$  represents a stiffness matrix from a discretization of the Laplace operator in $K_h$. The rectangular matrix $\mathbb{M}$ represents the $L^2-$ scalar product of functions from spaces $K_h$ and $\Lambda_h$. It holds
$$  \intO   \nabla v_h \cdot\nabla v_h  \, dx  = \textbf{v}^T \mathbb{K} \textbf{v}, \qquad  \intO  v_h \cdot \mu  \, dx  = \textbf{v}^T \mathbb{M} \boldsymbol{\mu} $$
for all $v_h \in K_h, \mu \in \Lambda_h$ and corresponding collumn vectors $\textbf{v}  \in \mathbb{R}^{|\NN|},  \boldsymbol{\mu} \in \mathbb{R}^{|\E|}$.
If we order nodes $\NN$ in a way that internal nodes $\NN_I$  precede Dirichlet nodes $\NN_D$ (no Neumann nodes are assumed for simplicity), we have the decomposition
$$\textbf{v}=(\textbf{v}_I, \textbf{v}_D)  \in \mathbb{R}^{|\NN_I|} \times \mathbb{R}^{|\NN_D|} $$
in Dirichlet and internals components and $|\NN|= |\NN_I|+ |\NN_D|$. Then \eqref{discrete_J_mu} can be rewritten  as
\begin{equation}\label{discrete_J_mu_components}
J_{\boldsymbol{\mu}}(\textbf{v})= \frac{1}{2} \begin{pmatrix} \textbf{v}_I   \\ \textbf{v}_D \end{pmatrix}^T
 \begin{pmatrix} \mathbb{K}_{I,I} &\mathbb{K}_{I,D}^T \\ \mathbb{K}_{I,D} & \mathbb{K}_{D,D} \end{pmatrix}
\begin{pmatrix} \textbf{v}_I  \\ \textbf{v}_D \end{pmatrix}
 +  \begin{pmatrix} \textbf{v}_I   \\ \textbf{v}_D \end{pmatrix}^T
 \begin{pmatrix} \mathbb{M}_{I} \\ \mathbb{M}_D \end{pmatrix}
\boldsymbol{\mu}.
\end{equation}
Note that the rectangular matrix $ \mathbb{K}_{ID}$ is the restriction of $ \mathbb{K}$ to its subblock with rows $\NN_I$ and columns  $\NN_D$. Therefore
$ \mathbb{K}_{I,I}$ and $ \mathbb{K}_{D,D}$ are not diagonal matrices. The rectangular matrices $\mathbb{M}_{I}$ and $\mathbb{M}_{D}$ are then restrictions of $\mathbb{M}$ to subblocks with rows  $\NN_I$ and  $\NN_D$ with all columns left.
The value of $\textbf{v}_D$ is known and given by Dirichlet boundary conditions. The direct computation
%is  computable from  $I(\boldsymbol{\mu})= \min_{\textbf{v}_I} J_{\boldsymbol{\mu}}(\textbf{v})$
reveals
\begin{equation}\label{dual_discrete_I_evaluation_form}
I^*(\boldsymbol{\mu})= \frac{1}{2} \textbf{v}_D^T \, \mathbb{K}_{D,D} \, \textbf{v}_D +  \textbf{v}_D^T \, \mathbb{M}_{D} \, \boldsymbol{\mu}
- \frac{1}{2} ( \mathbb{K}_{I,D} \, \textbf{v}_D   +  \mathbb{M}_{I} \, \boldsymbol{\mu}  )  \mathbb{K}_{I,I}^{-1}  ( \mathbb{K}_{I,D} \, \textbf{v}_D   +  \mathbb{M}_{I} \, \boldsymbol{\mu}  ).
\end{equation}
In addition to it, for a given $\boldsymbol{\mu}$, the component $\textbf{v}_I$ minimizing the functional \eqref{discrete_J_mu_components} satisfies
\begin{equation}\label{formula_vi_from_mu}
  \textbf{v}_I =-\mathbb{K}_{I,I}^{-1}  ( \mathbb{K}_{I,D} \, \textbf{v}_D   +  \mathbb{M}_{I} \, \boldsymbol{\mu}) \qquad \mbox{or equivalently}  \quad -\mathbb{K}_{I,I} \textbf{v}_I = ( \mathbb{K}_{I,D} \, \textbf{v}_D   +  \mathbb{M}_{I} \, \boldsymbol{\mu}).
\end{equation}

This formula is applied for the reconstruction of $\textbf{v}_I$ from $ \boldsymbol{\mu}$.
 Since $\mathbb{K}_{I,I}$ is  a sparse matrix, its inverse  $\mathbb{K}_{I,I}^{-1}$ is  a dense  matrix. Then, the right part of  of \eqref{formula_vi_from_mu} is exploited in practical evolutions including \eqref{dual_discrete_I_evaluation_form}. The matrix $\mathbb{K}_{I,I}^{-1}$ is positive definite as well as $\mathbb{K}_{I,I}$ and the functional $I^*(\cdot)$ is concave and contains  quadratic and linear terms only. Thus the functional $-I^*(\cdot)$ is convex and its minimum $\lambda_h$ from \eqref{dual_discrete} is represented by a column vector $\boldsymbol{\lambda} \in  \mathbb{R}^{|\T|} $ and solves a quadratic programming (QP) problem with box constraints
$$  -I^*(\boldsymbol{\lambda}) = \min -I^*(\boldsymbol{\mu}), \qquad \mbox{where } \{ \boldsymbol{\mu} \}_j \in [-\alpha_{-}, \alpha_{+}] \mbox{ for all } j \in \{1,\dots,  |\T|  \}.$$
The corresponding solution $u_{\lambda_h} \in K_h$ is then represented by a collumn vector
$\textbf{u}_{\boldsymbol{\lambda}}=(\textbf{u}_{{\boldsymbol{\lambda}}_I}, \textbf{u}_{{\boldsymbol{\lambda}}_D}) \in  \mathbb{R}^{|\NN|} $ and $\textbf{u}_{{\boldsymbol{\lambda}}_I}$ solves \eqref{formula_vi_from_mu} for $\textbf{v}=\boldsymbol{\lambda}$. The solutions steps above are summarized in Algorithm \ref{alg:dual}.
%Since QP is an iterative method, it is terminated with respect to its stopping criteria and provide only an approximation $\mu_h$ to  $\lambda_h$

\begin{algorithm}
\noindent Let discretization matrices $ \mathbb{K}  \in \mathbb{R}^{|\NN| \times |\NN|} $ and $\mathbb{M} \in \mathbb{R}^{|\NN| \times |\T|} $ be given with their subblocks $ \mathbb{K}_{I,I} \in \mathbb{R}^{|\NN_I| \times |\NN_I|}  , \mathbb{K}_{D,D}\in \mathbb{R}^{|\NN_D| \times |\NN_D|} , \mathbb{K}_{I,D} \in \mathbb{R}^{|\NN_I| \times |\NN_D|} $ and $ \mathbb{M}_{I} \in \mathbb{R}^{|\NN_I| \times |\T| } , \mathbb{M}_{D} \in \mathbb{R}^{|\NN_D| \times |\T| }$. Let  $\textbf{v}_D \in  \mathbb{R}^{|\NN_D|}$ be vector of prescribed Dirichlet values in nodes $\NN_D$. Then:
\begin{itemize}
\item[(i)] find the vector of Lagrange multipliers from the quadratic minimization problem
$$\boldsymbol{\lambda} = \argmin_{ \boldsymbol{\mu} \in \mathbb{R}^{|\T|}  }
\left(-  \frac{1}{2} \textbf{v}_D^T \, \mathbb{K}_{D,D} \, \textbf{v}_D -  \textbf{v}_D^T \, \mathbb{M}_{D} \, \boldsymbol{\mu}
+ \frac{1}{2} ( \mathbb{K}_{I,D} \, \textbf{v}_D   +  \mathbb{M}_{I} \, \boldsymbol{\mu}  )  \mathbb{K}_{I,I}^{-1}  ( \mathbb{K}_{I,D} \, \textbf{v}_D   +  \mathbb{M}_{I} \, \boldsymbol{\mu})   \right) $$
under box constraints $ \{ \boldsymbol{\mu} \}_j \in [-\alpha_{-}, \alpha_{+}] \mbox{ for all } j \in \{1,\dots,  |\T|  \},$

\item[(ii)] reconstruct the solution vector $\textbf{v}_I  \in \mathbb{R}^{|\NN_I|} $ from
$\quad -\mathbb{K}_{I,I} \textbf{v}_I = ( \mathbb{K}_{I,D} \, \textbf{v}_D   +  \mathbb{M}_{I} \, \boldsymbol{\mu})$.

\item[(iii)] output $\lambda_h$ and $u_{\lambda_h}$ represented by vectors $\boldsymbol{\lambda}$ and $\textbf{v}=(\textbf{v}_I, \textbf{v}_D)$. 
\end{itemize}
\caption{Quadratic programing for Langrange multipliers.} \label{alg:dual}
\end{algorithm}

\subsection{Minimization of the functional error majorant}
For a given approximation $u_{\lambda_h} \in   K_h $,   the majorant value $\mathcal{M_+}(u_{\lambda_h} ; \beta, \eta^{*}, \mu)$ majorizes the value $J(u_{\lambda_h}) - J(u)$. The majorant $\mathcal{M_+}(u_{\lambda_h} ; \beta, \eta^{*}, \mu)$ can be minimized with respect to its free arguments $\beta>0, \eta^{*} \in Y^*, \mu \in \Lambda$ in order to obtain the sharp upper bound.
The fields $\eta^{*} \in Y^*, \mu \in \Lambda$ can be sought on a mesh $\T_{\tilde h}$ with a different mesh size $\tilde h$. Choosing very small mesh size $\tilde{h} \ll h$ leads to sharper bounds but higher computational costs. Here we consider the same mesh size $\tilde{h}=h$ for simplicity,  $$\eta^{*} \in Y^*_h, \quad \mu \in \Lambda_h.$$
We use the successive minimization algorithm described in Algorithm \ref{alg:majorant_minimzation}.
\begin{algorithm}[h]
\noindent Let $k=0$ and let initial $\beta_{0}>0$ and $\mu_{0} \in \Lambda_h$ be given. Then:
\begin{itemize}
\item[(i)] find an iteration $\eta^{*}_{k+1}    \in Y^*_h$ such that
$
\eta^{*}_{k+1}=\argmin\limits_{\eta^{*}\in Y^*_h} \mathcal{M_+}(u_{\lambda_h};\beta_k,\mu_k,\eta^{*}),
$
\item[(ii)] find $\mu_{k+1}\in \Lam_h$ such that
$
\mu_{k+1}=\argmin\limits_{\mu\in \Lambda_h} \mathcal{M_+}(u_{\lambda_h};\beta_k,\mu,\eta^{*}_{k+1}),
$
\item[(iii)] find $\beta_{k+1}>0$ such that
$
\beta_{k+1}=\argmin\limits_{\beta \in \R_+ } \mathcal{M_+}(u_{\lambda_h};\beta,\mu_{k+1},\eta^{*}_{k+1}),
$
\item[(iv)] set $k:=k+1$ are repeat (i)\,--\,(iii) until convergence. Then,  output $\eta_h^{*}:=\eta^{*}_{k+1}$ and $\mu_h:=\mu_{k+1}$.
\end{itemize}
\caption{Majorant minimization algorithm.} \label{alg:majorant_minimzation}
\end{algorithm}

The step (i)   corresponds to the solution of a linear system of equations
\begin{equation}\label{eq:system_y}
\left[(1+\beta_k){\mathbb M}^{RT0}+C^2_{\Omega}\left(1+\frac{1}{\beta_k}\right){\mathbb K}^{RT0}\right]{\boldsymbol{\eta}^*_{k+1} }=(1+\beta_k)\textbf{c} +C^2_{\Omega}\left(1+\frac{1}{\beta_k}\right)\textbf{d}
\end{equation}
for a column vector $\boldsymbol{\eta}^*_{k+1}  \in \mathbb{R}^{|\E|} $. Here, ${\mathbb K}^{RT0}, {\mathbb M}^{RT0} \in \mathbb{R}^{|\E|\times |\E|} $ are stiffness and mass matrices corresponding to $RT_0(\mathcal T_h)$ elements, described together with vectors $\textbf{c}, \textbf{d}$ in Section 4 of \cite{HaVa}. The minimal argument $\mu_{k+1} \in \Lambda_{h}$ in step (ii)  is  locally computed on every triangle $T \in \T_h$ from the formula
\begin{equation}\label{mu_upgrade}
\mu_{k+ 1} |_T
=
\mathcal{P}_{[ -\alpha_{-}, \alpha_{+}  ]}
\left(
{\rm{div}}\,\eta^{*}_{k+1}|_T
+\frac{\overline{u}_{\lambda_h}|_T}
{C_{\Omega}^{2}\left(1+\frac{1}{\beta_{k}}\right)}
\right) ,
\end{equation}
where $\mathcal{P}_{[ -\alpha_{-}, \alpha_{+}  ]}$ is the projection on the convex set $[ -\alpha_{-}, \alpha_{+}  ]$ and $\overline{u}_{\lambda_h}|_T$ means an averaged value of $u_{\lambda_h}$  over a triangular element $T$. The minimization in step (iii) leads to the explicit relation
\begin{equation}  \label{beta_upgrade}
\beta_{k+1}=\frac{\|{\rm{div}} \,\eta^{*}_{k+1} - \mu_{k+1}\|_\Omega}{ \| \nabla v-\eta^{*}_{k+1} \|_{\Omega}} \,.
\end{equation}
\begin{remark}[Choice of initial $\mu_0$]
We recall the approximation $\lambda_h$ is taken as  an initial approximation $\mu_0$, which can speed up the convergence significantly \cite{HaVa}.
\end{remark} 

\section{Numerical examples} \label{sec:numerics}
In this section we elaborate two numerical examples, i.e. Example I, and Example II, with known and unknown exact solutions.

\subsection{Example I with known exact solution} \label{subsec:ex1}
This example is introduced in        \cite{Bo}; it is  also tested for  one dimensional case in   \cite{ReVa2015}. Here, we consider it  in two dimensional   and assume a rectangular domain
\begin{equation}  \label{example1_domain}
\Omega=X \times Y:=(-1, 1) \times (0, 1),
\end{equation}
and constant coefficients
\begin{equation}  \label{example1_coefficients}
\alpha_{-}=\alpha_{+}=8.
\end{equation}
The two phase obstacle problem   \eqref{H}  is supplied with the  Dirichlet boundary conditions
\begin{equation}  \label{example1_BC}
u(-1, y)=-1, \quad u(1, y)=1 \qquad \forall y \in Y
\end{equation}
and homogeneous Neumann boundary conditions
\begin{equation}  \label{example1_BC_Neumann}
\frac{\partial u}{\partial x}(x,0)=\frac{\partial u}{\partial x}(x,1)=0 \qquad \forall x \in X.
\end{equation}

\begin{figure}
\center
\begin{minipage}{0.45\textwidth}
\includegraphics[width=\textwidth]{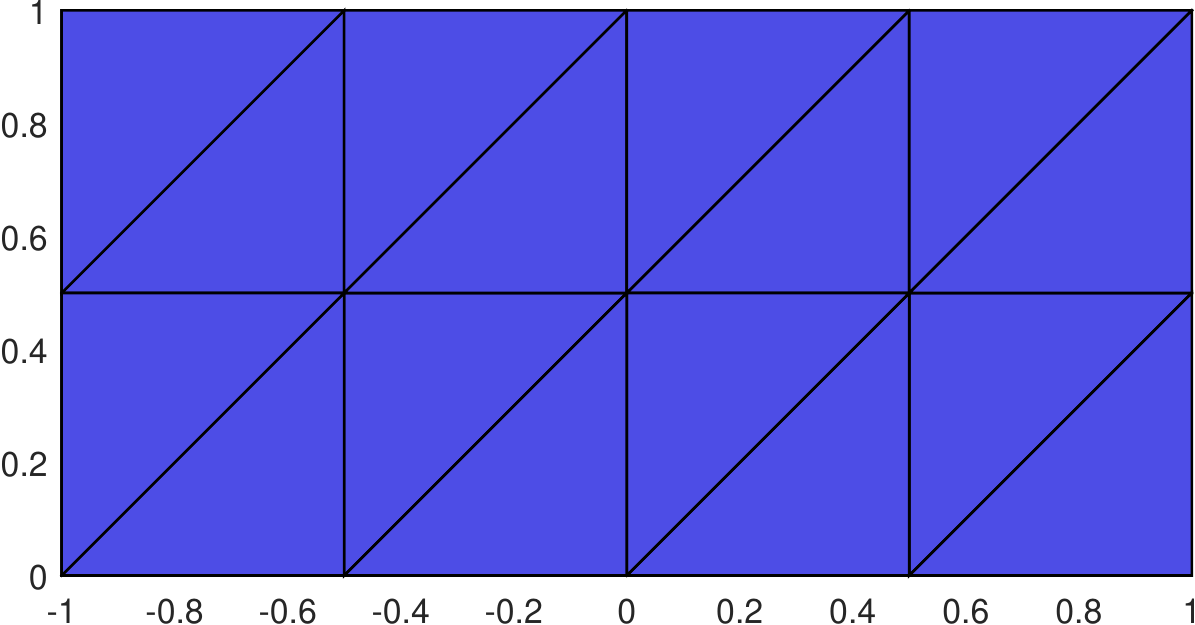}
\end{minipage}
\hspace{1cm}
\begin{minipage}{0.45\textwidth}
\includegraphics[width=\textwidth]{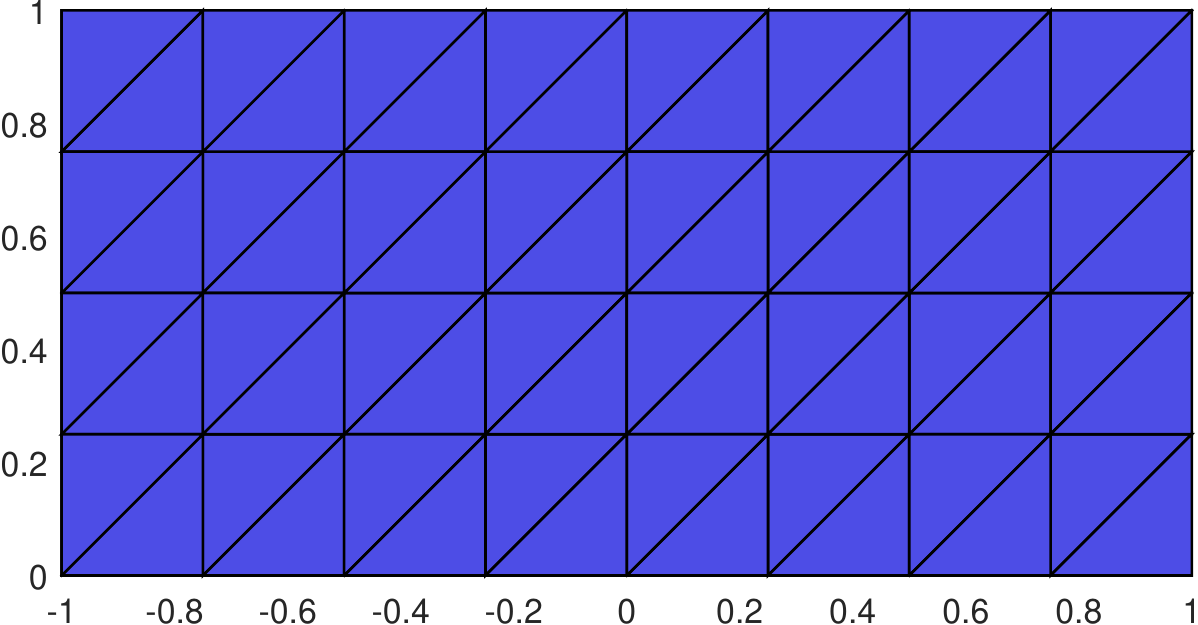}
\end{minipage}
\caption{Example I - level 1 and level 2 nested triangular meshes. Note that there are triangular node lying on the exact free boundary is given by lines $x=\pm 0.5$.}
\label{figure_meshes_exampleI}
\vspace{1cm}
\begin{minipage}{0.45\textwidth}
\includegraphics[width=\textwidth]{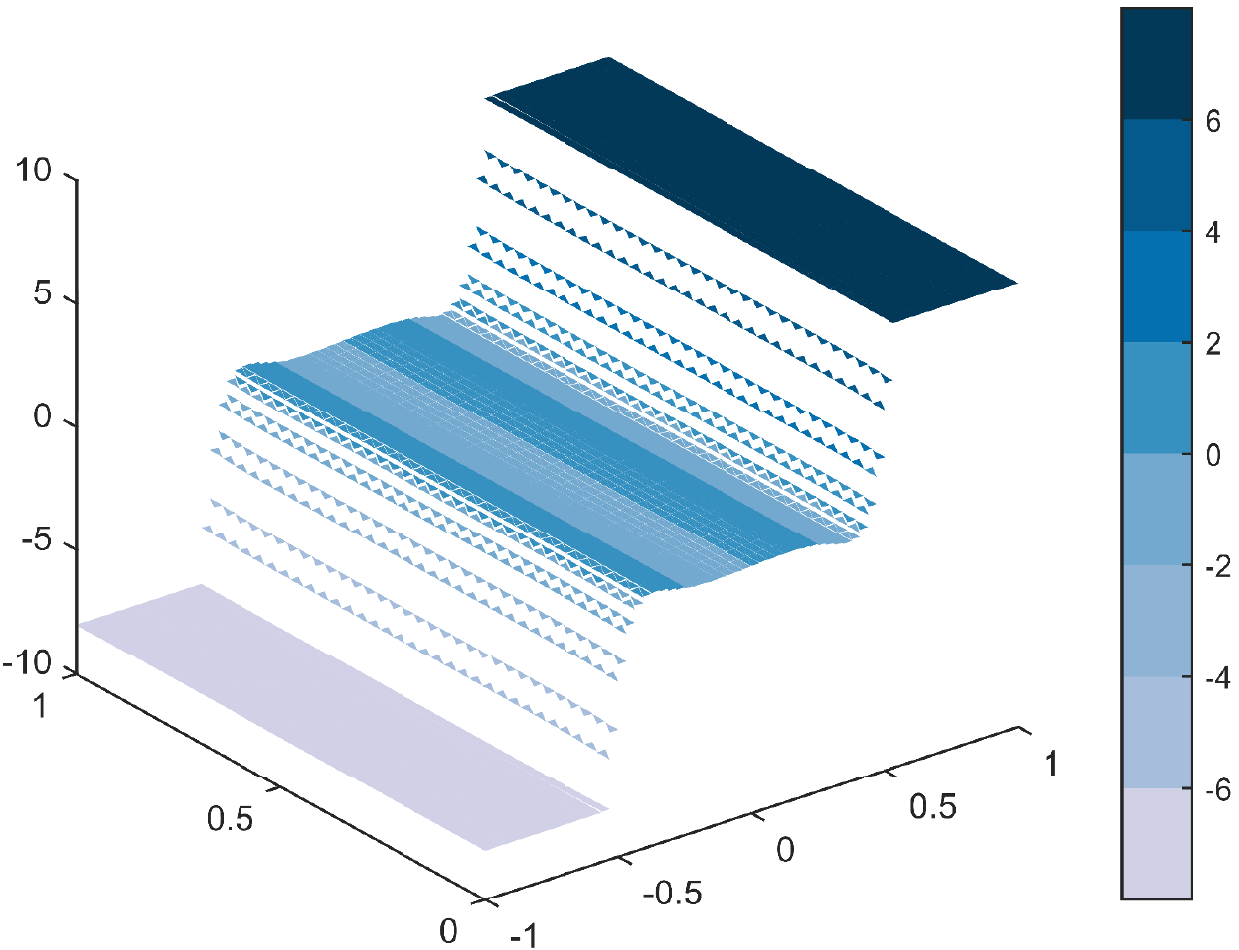}
\end{minipage}
\hspace{1cm}
\begin{minipage}{0.45\textwidth}
\includegraphics[width=\textwidth]{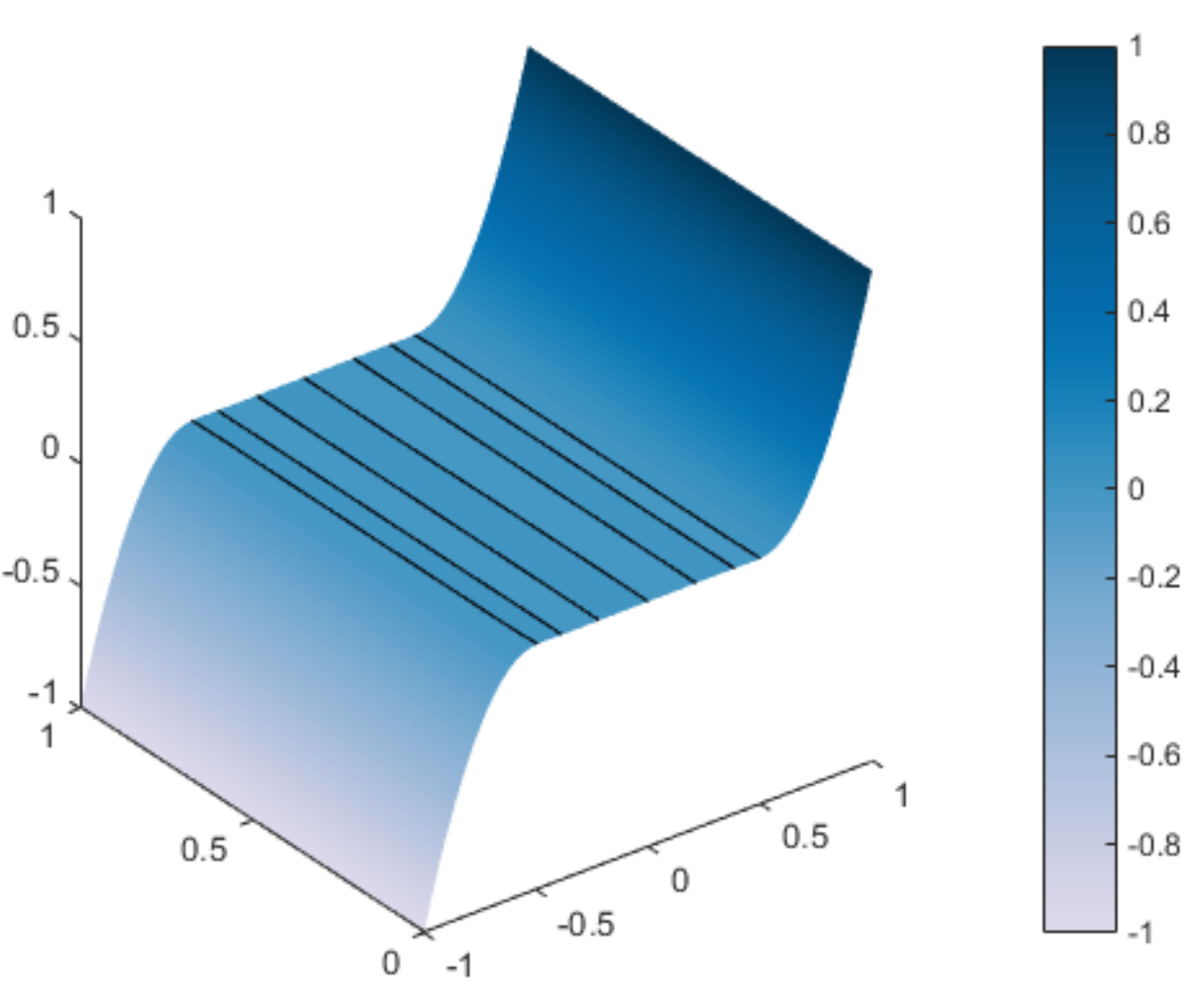}
\end{minipage}
\caption{Example I - approximations: multiplier $\lambda_h \in \Lambda_h$ (left) and the corresponding solution $u_{\lambda_h} \in K_h$ (right) computed on level 5 triangular mesh (referred to as level 5 in Table \ref{table_exampleI}). The multiplier approximation $\lambda_h$  indicates an approximative free boundary, the exact free boundary is given by lines $x=\pm 0.5$. Full contour lines of $u_{\lambda_h}$ at values  $\pm 0.0001$ are additionally displayed (right).
%{\color{red}  I see there are more than   $\pm 0.0001$ ? }
}
\label{figure_solution_exampleI}
\vspace{1cm}
\begin{minipage}{0.32\textwidth}
\includegraphics[width=\textwidth]{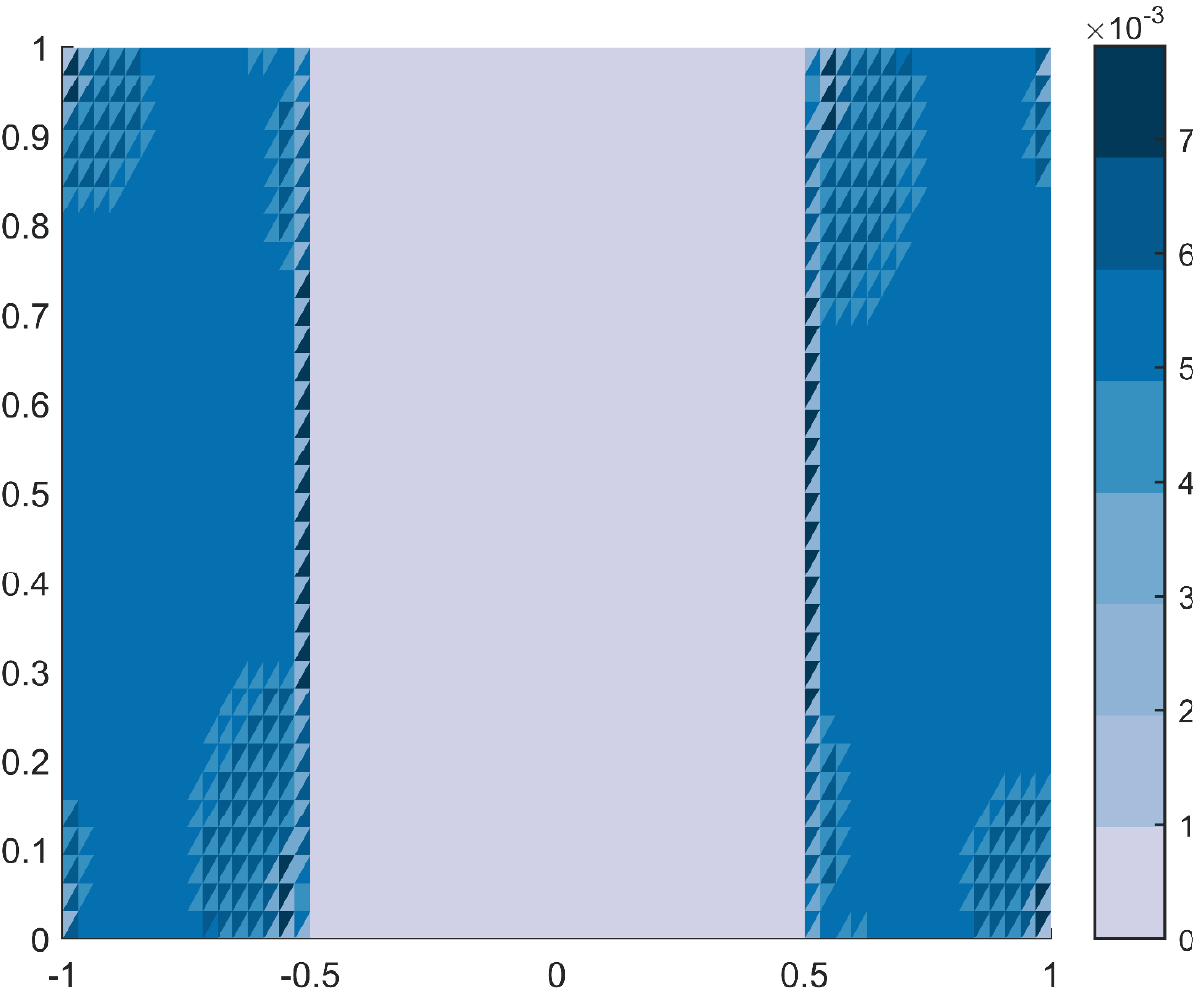}
\end{minipage}
\begin{minipage}{0.32\textwidth}
\includegraphics[width=\textwidth]{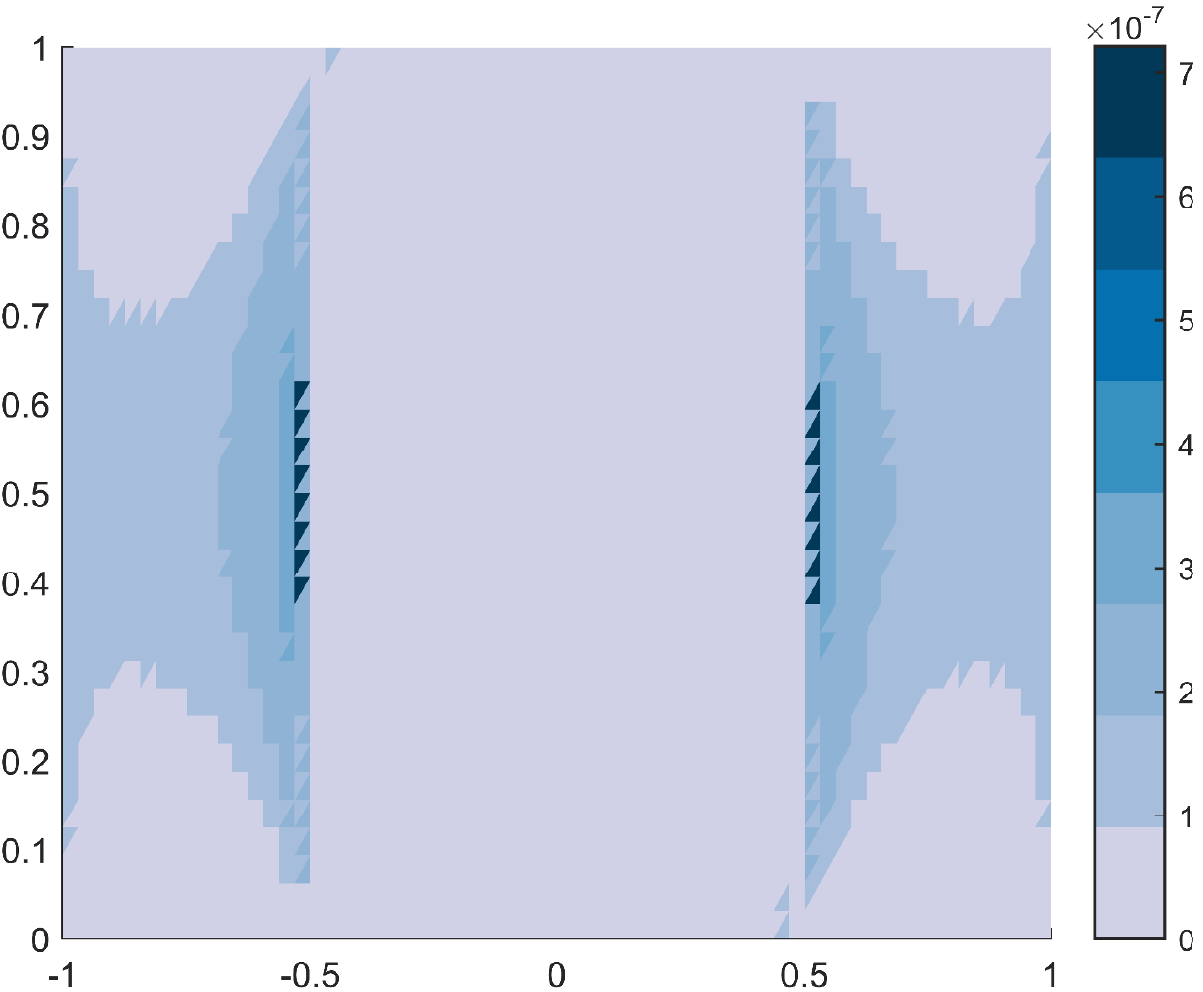}
\end{minipage}
\begin{minipage}{0.32\textwidth}
\includegraphics[width=\textwidth]{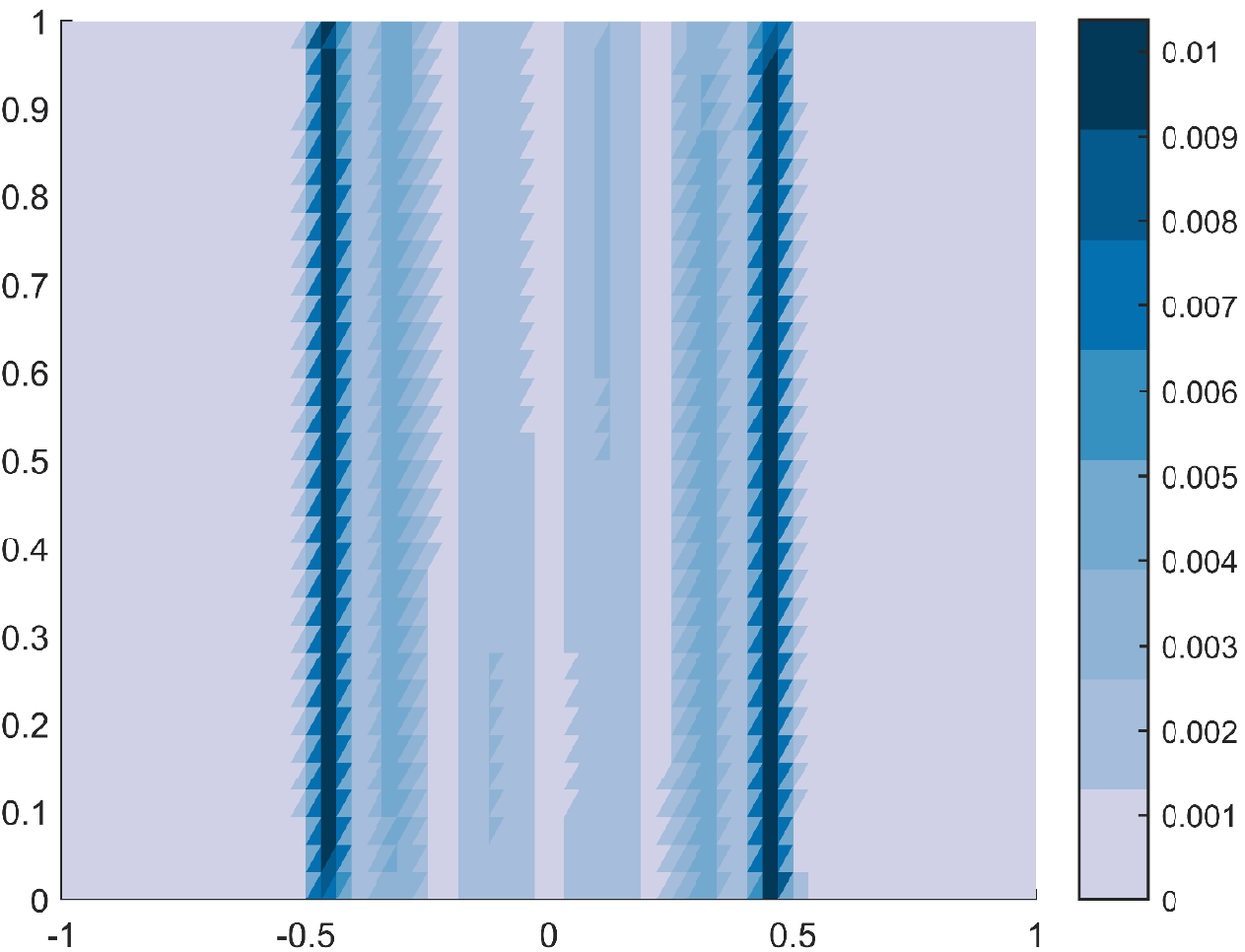}
\end{minipage}
\caption{Example I - distribution of the majorant parts  $\mathcal{M_+}_1$ (left), $\mathcal{M_+}_2$ (middle), $\mathcal{M_+}_3$ (right)
computed on level 5 triangular mesh (referred to as level 5 in Table \ref{table_exampleI}).}
\label{figure_majorantParts_exampleI}
%\vspace{1cm}
%\begin{minipage}{0.45\textwidth}
%\includegraphics[width=\textwidth]{example_I_error}
%\end{minipage}
%\hspace{1cm}
%\begin{minipage}{0.45\textwidth}
%\includegraphics[width=\textwidth]{example_I_majorant}
%\end{minipage}
%\caption{Example I - distribution of the error $\frac{1}{2}\|\nabla(u_{\lambda_h}-u)\|^{2}_{\Omega,A}$ (left) and of the majorant $\mathcal{M_+}$  (right) computed on a triangular mesh (referred to as level 5 mesh in Table \ref{table_exampleI}).}
%\label{figure_error_and_majorant_exampleI}
\end{figure}

The exact solution $u \in K$ is given by the relation independent of $y \in Y$,
\begin{equation}\label{u_exact_ExampleI}
u(x,y)=\left\{\begin{array}{ll} -4x^2-4x-1 , \quad & x \in X_{-}:= [-1, -0.5], \\
                                 0, \quad & x \in X_0:= [-0.5, 0.5], \\
                                 4x^2-4x+1, \quad & x \in X_{+}:= [0.5, 1]
                                    \end{array} \right.
\end{equation}
and its (exact) energy is $J(u)=5\frac{1}{3}.$
The (exact) free boundary is characterized by two lines
$$ (\pm 0.5, y), \qquad \mbox{ where } y \in Y.$$
The (exact) Lagrange multiplier $\lambda \in \Lam $ is then given by
\begin{equation}\label{lambda_exact_ExampleI}
 \lambda(x,y)=\left\{\begin{array}{ll} -\alpha_{-} , \quad & x \in X_{-}, \\
                                 0, \quad & x \in X_{0}, \\
                                 \alpha_{+}, \quad & x \in X_{+}
                                    \end{array} \right.
\end{equation}
and it is a discontinuous function with a jump on the free boundary. We compute approximation pairs $(\lambda_h,u_{\lambda_h}) \in \Lam_h \times K_h$ for a sequence of nested uniformly refined meshes. Levels 1 and 2 meshes are depicted in Figure \ref{figure_meshes_exampleI}. Since some dicretization nodes are lying exactly on the free boundary, there might be a chance to reconstruct the free boundary exactly from approximative solutions. A finer (level 5) approximation pair $(\lambda_h,u_{\lambda_h}) \in \Lam_h \times K_h$ computed from the dual-based solver is shown in Figure \ref{figure_solution_exampleI}.  The approximative Lagrange multiplier field $\lambda_h$ however only approximates the exact free boundary. 
%Our approximation $u_{\lambda_h}$ is computed as a $P_1(\mathcal T_h)$ function but $u$ is known to be piecewise quadratic {\color{red} piecewise q%uadratic  $u$ is $C^{1,1}$ so I don't know piecewise quadratic }. It is expected that if $P_2(\mathcal T_h)$ elements were used for $u_{\lambda_h}$ and %$P_0(\mathcal T_h)$ elements for $\lambda_h$, the exact free boundary and the exact energy could be computed{\color{red}Sounds like that you are not %sure, or you have doubt! If that is a fact, then remove this part}. 
To the given approximation pair $(\lambda_h,u_{\lambda_h})$, a functional majorant is optimized using 10000 iterations of Algorithm 1 ( we set $\mu_0=\lambda_h$). To get more insight on the majorant behaviour, we display space densities of all three additive majorant subparts
$$\mathcal{M_+}_1(u_{\lambda_h}; \beta, \eta^{*}), \qquad \mathcal{M_+}_2(\beta, \eta^{*}, \mu   ), \qquad \mathcal{M_+}_3(u_{\lambda_h}; \mu)$$ separately in Figure \ref{figure_majorantParts_exampleI}. The amplitudes of $\mathcal{M_+}_2$  are significantly lower than amplitudes of $\mathcal{M_+}_1$ and $\mathcal{M_+}_3$, but the value of $\mathcal{M_+}_3$ is still relatively high. The high value of $\mathcal{M_+}_3$ indicates that the exact free boundary is not sufficiently resolved yet and the density of $\mathcal{M_+}_3$ seem to be a reasonable indicator of the exact free boundary.

%\begin{remark}[convergence of majorant parts]
%Theorem 3.1 of \cite[ReVa2015] state that majorant for the case of the optimized exact parameters $\eta=\nabla u$, %$\mu=\lambda$ is exactly equal to
%\end{remark}

\begin{table}
\begin{center}
\begin{small}
\begin{tabular}{|c|c|c|c|c|c|c|c|c|}
\hline
$\mbox{level}$  & $|\NN|$  & $J(u_{\lambda_h}) \, (I^*(\lambda_h))$ & $J(u_{\lambda_h}) -J(u)$
 & $\mathcal{M_+}(u_{\lambda_h},\cdot)$   & $\mathcal{M_+}_1(\cdot)$  & $\mathcal{M_+}_2(\cdot)$  & $\mathcal{M_+}_3(\cdot)$
%& $\frac{\sqrt{2 \mathcal{M}(v,\dots)}}{\|\nabla(v-u)\|_{\Omega,A}}$
\\
%   &  &   &  &  efficiency \\
\hline
%0 & 6 & 8.9965 (8.4074) & 3.66e+00 & 1.05e+01 & 1.53e+00 & 1.31e-01 & 8.88e+00 \\
1 & 15 & 6.0383 (5.9975) & 7.05e-01 & 1.74e+00 & 1.33e+00 & 3.00e-02 & 3.80e-01 \\
2 & 45 & 5.5030 (5.5000) & 1.70e-01 & 3.79e-01 & 3.47e-01 & 9.44e-04 & 3.15e-02 \\
3 & 153 & 5.3765 (5.3750) & 4.32e-02 & 9.14e-02 & 8.41e-02 & 2.47e-05 & 7.28e-03 \\
4 & 561 & 5.3481 (5.3435) & 1.47e-02 & 2.76e-02 & 2.10e-02 & 6.08e-07 & 6.50e-03 \\
5 & 2145 & 5.3390 (5.3355) & 5.65e-03 & 8.62e-03 & 5.21e-03 & 7.13e-08 & 3.41e-03 \\
\hline
\end{tabular}
\caption{Computations of Example I on various uniform triangular meshes. }\label{table_exampleI}
\end{small}
\end{center}
\end{table}
Computations on all nested uniformly refined triangular meshes are summarized in  \ref{table_exampleI}. The dual energy $I^*(\lambda_h)$ and primal primal energy $J(u_{\lambda_h})$ converge to the exact energy $J(u)$ as $h \rightarrow 0$. Since we work with nested meshes, we additionally have
$$J(u_{\lambda_h}) \searrow J(u) \qquad \mbox{ or equivalently } \qquad J(u_{\lambda_h}) - J(u) \searrow 0 .$$
The difference of energies $ J(u_{\lambda_h}) - J(u) \searrow 0$ is bounded from above by the majorant value $\mathcal{M_+}(u_{\lambda_h},\dots)$  as stated in the majorant estimate \eqref{eq:majorant_estimate}.

\begin{remark}[Extension to mixed Dirichlet - Neumann boundary conditions]
This example assumes both Dirichlet and Neumann boundary conditions, but only Dirichlet boundary conditions are considered in $K$. The dual based solver for a double-phase problem can still be applied, with Neumann nodes  $\NN_N$ being added to internal nodes $\NN_I$. The majorant estimate \eqref{eq:majorant_estimate} is valid with the same majorant form \eqref{majorant}, but the flux  $\eta^{*}\in Q_h$ must satisfy an extra condition $\eta^{*} \cdot n=0$ on a Neumann boundary, where $n$ is a normal vector to the boundary. This condition means that components of   $\boldsymbol{\eta}^*_{k+1}$ from \eqref{eq:system_y} corresponding to Neumann edges $\E_N$ must be equal to zero.
\end{remark}

\subsection{Example II} \label{subsec:ex2}
\begin{figure}
\center
\begin{minipage}{0.4\textwidth}
\includegraphics[width=0.9\textwidth]{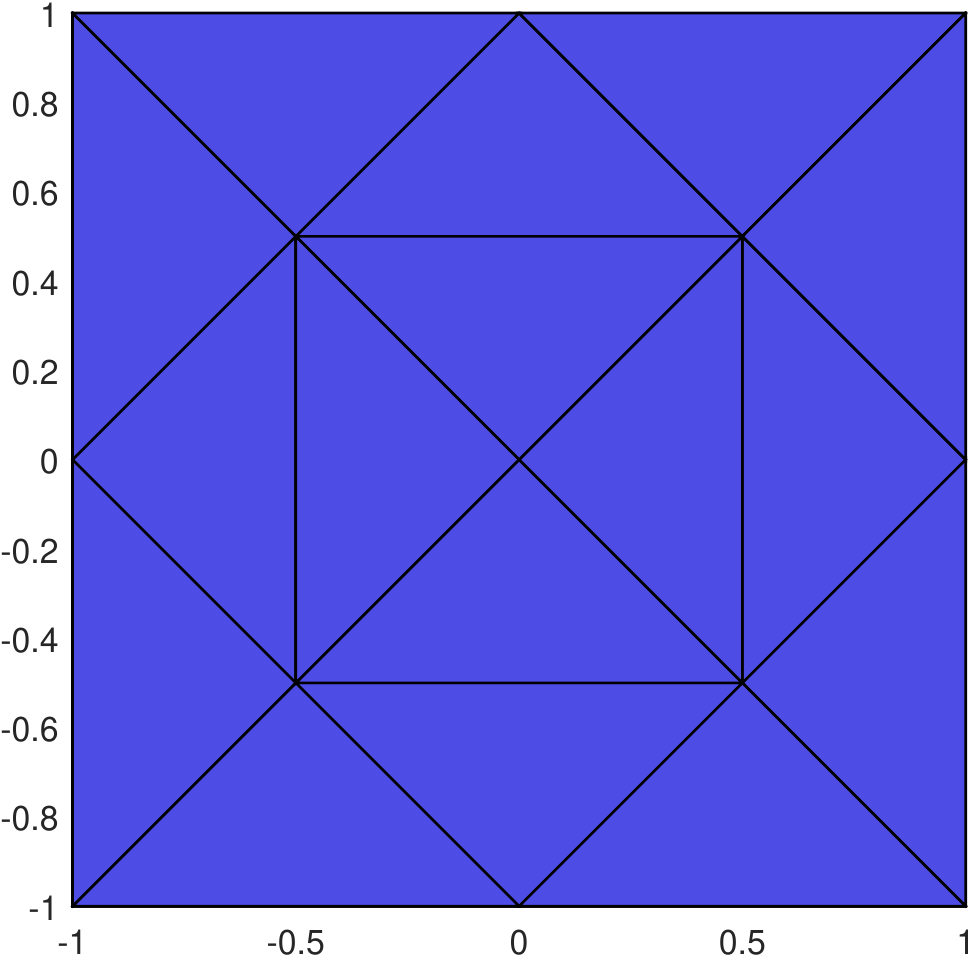}
\end{minipage}
\hspace{2cm}
\begin{minipage}{0.4\textwidth}
\includegraphics[width=0.9\textwidth]{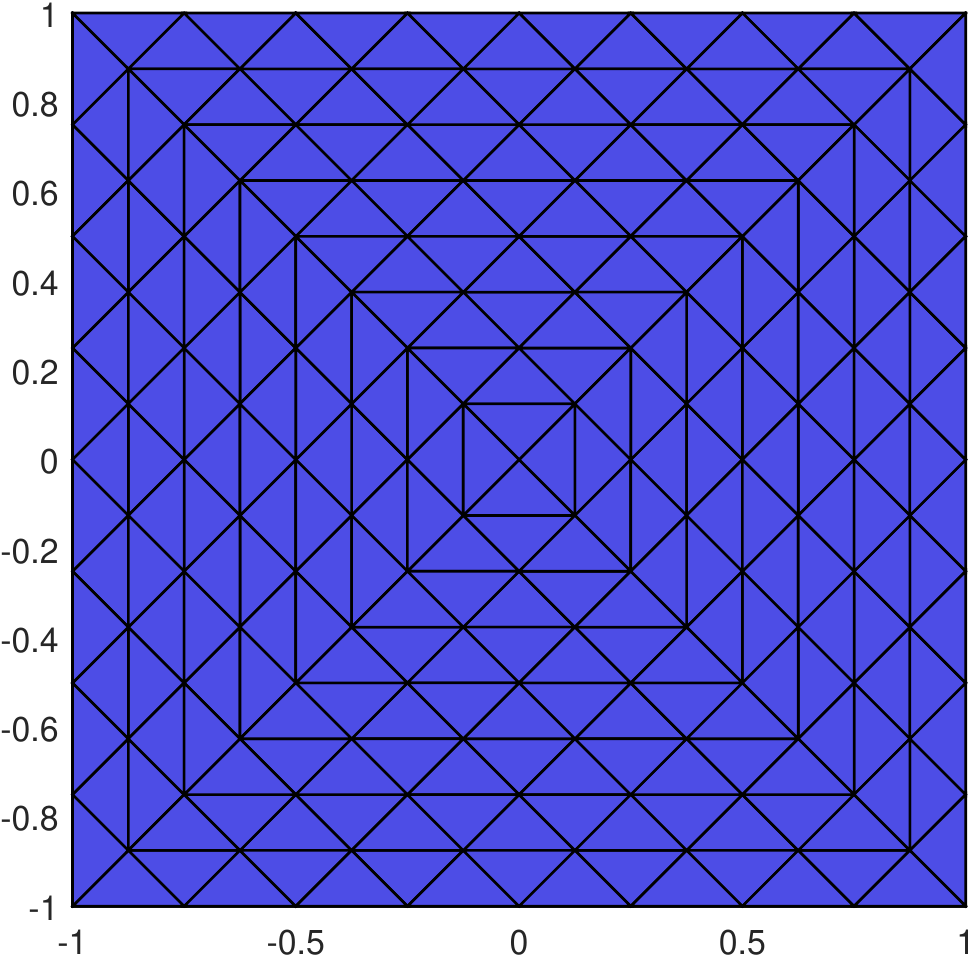}
\end{minipage}
\caption{Example II - level 1 and level 2 nested triangular meshes.}
\label{figure_meshes_exampleII}
\vspace{1cm}
\begin{minipage}{0.45\textwidth}
\includegraphics[width=\textwidth]{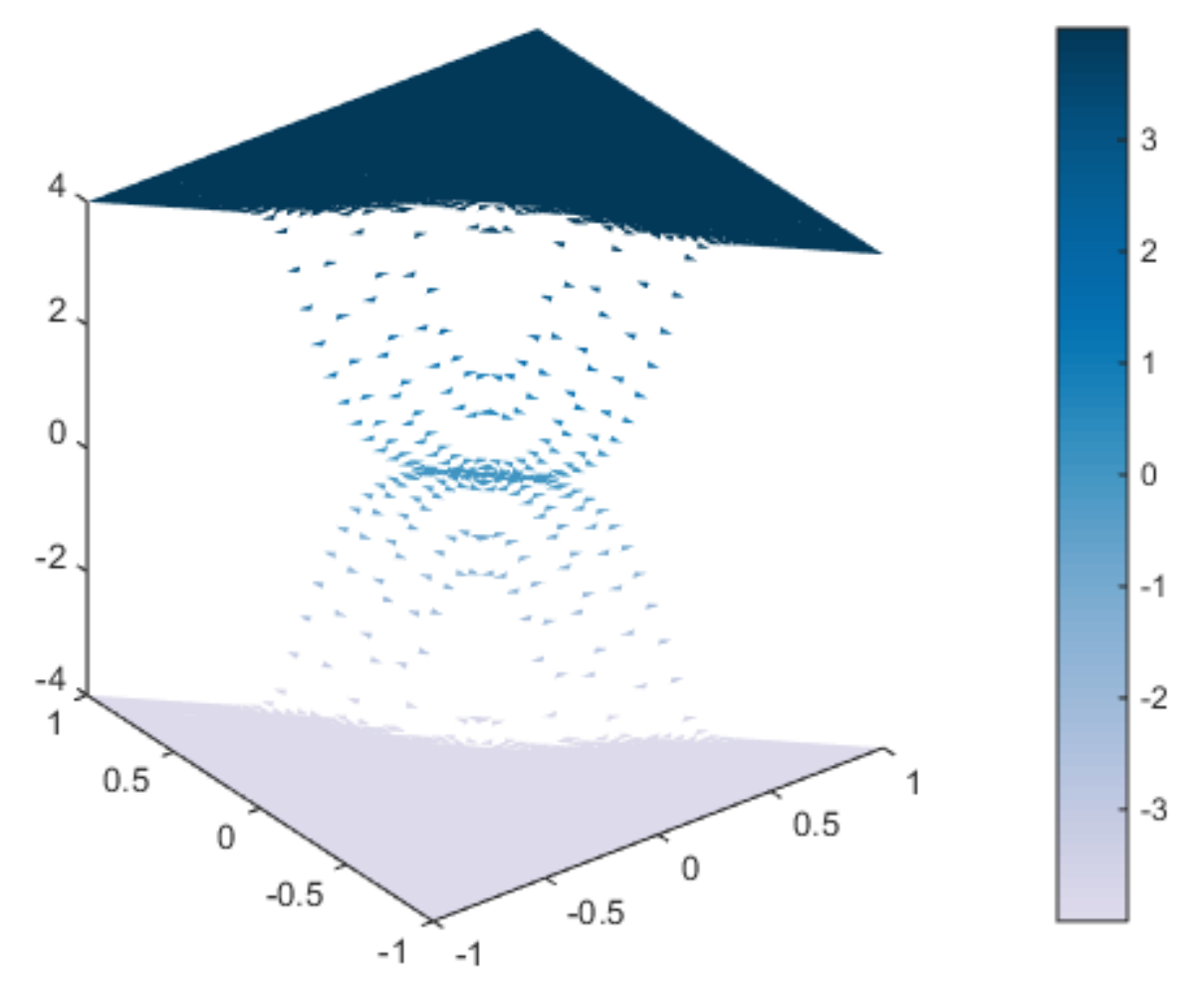}
\end{minipage}
\hspace{1cm}
\begin{minipage}{0.45\textwidth}
\includegraphics[width=\textwidth]{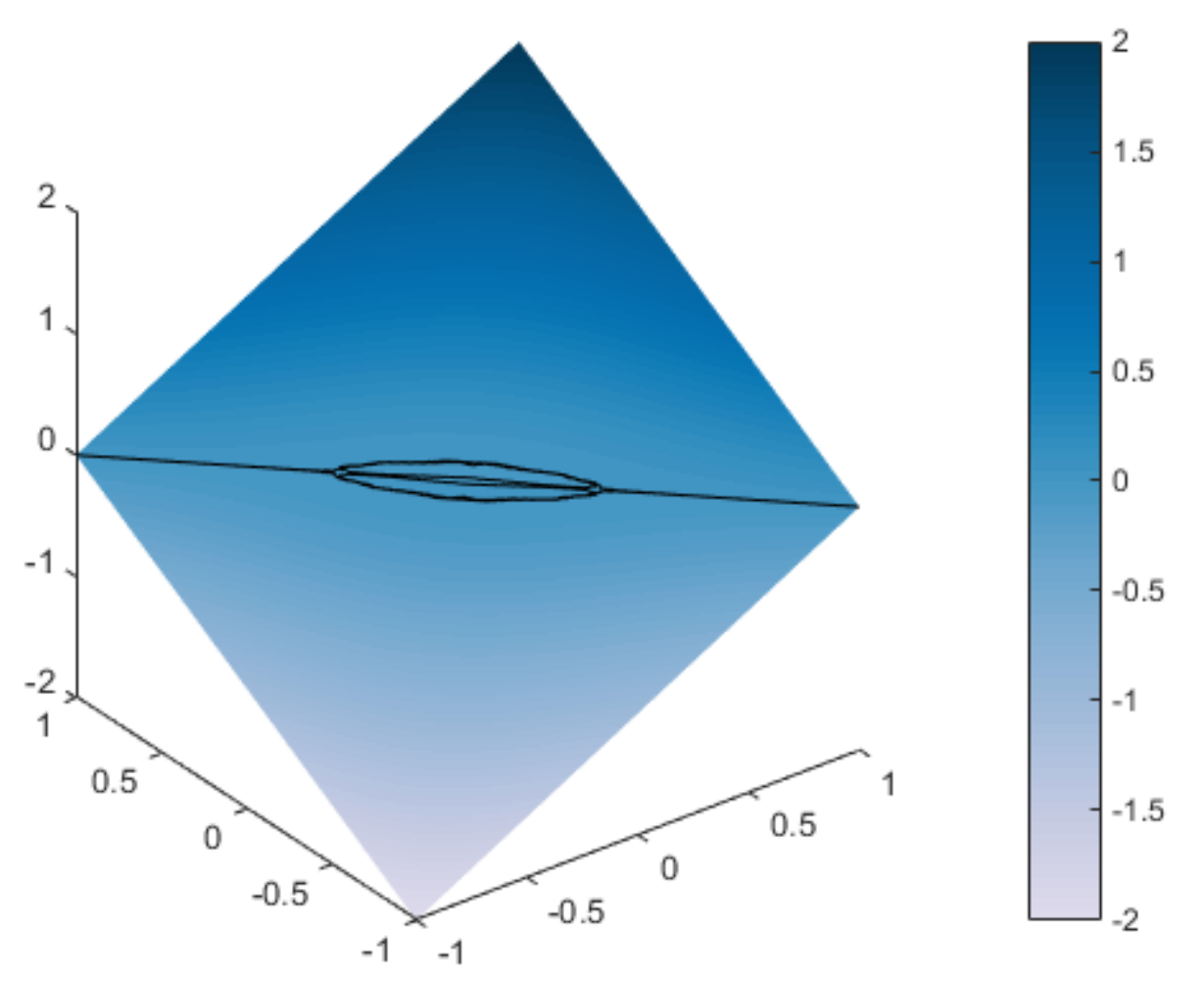}
\end{minipage}
\caption{Example II - approximations: multiplier $\lambda_h \in \Lambda_h$ (left) and the corresponding solution $u_{\lambda_h} \in K_h$ (right) computed on level 5 triangular mesh (referred to as level 5 in Table \ref{table_exampleII}). The multiplier approximation $\lambda_h$ (left) indicates an approximative free boundary, the exact free boundary  is unknown. Full contour lines of $u_{\lambda_h}$ at values  $\pm 0.0001$ are additionally displayed (right).
   }
\label{figure_solution_exampleII}
\vspace{1cm}
\begin{minipage}{0.32\textwidth}
\includegraphics[width=\textwidth]{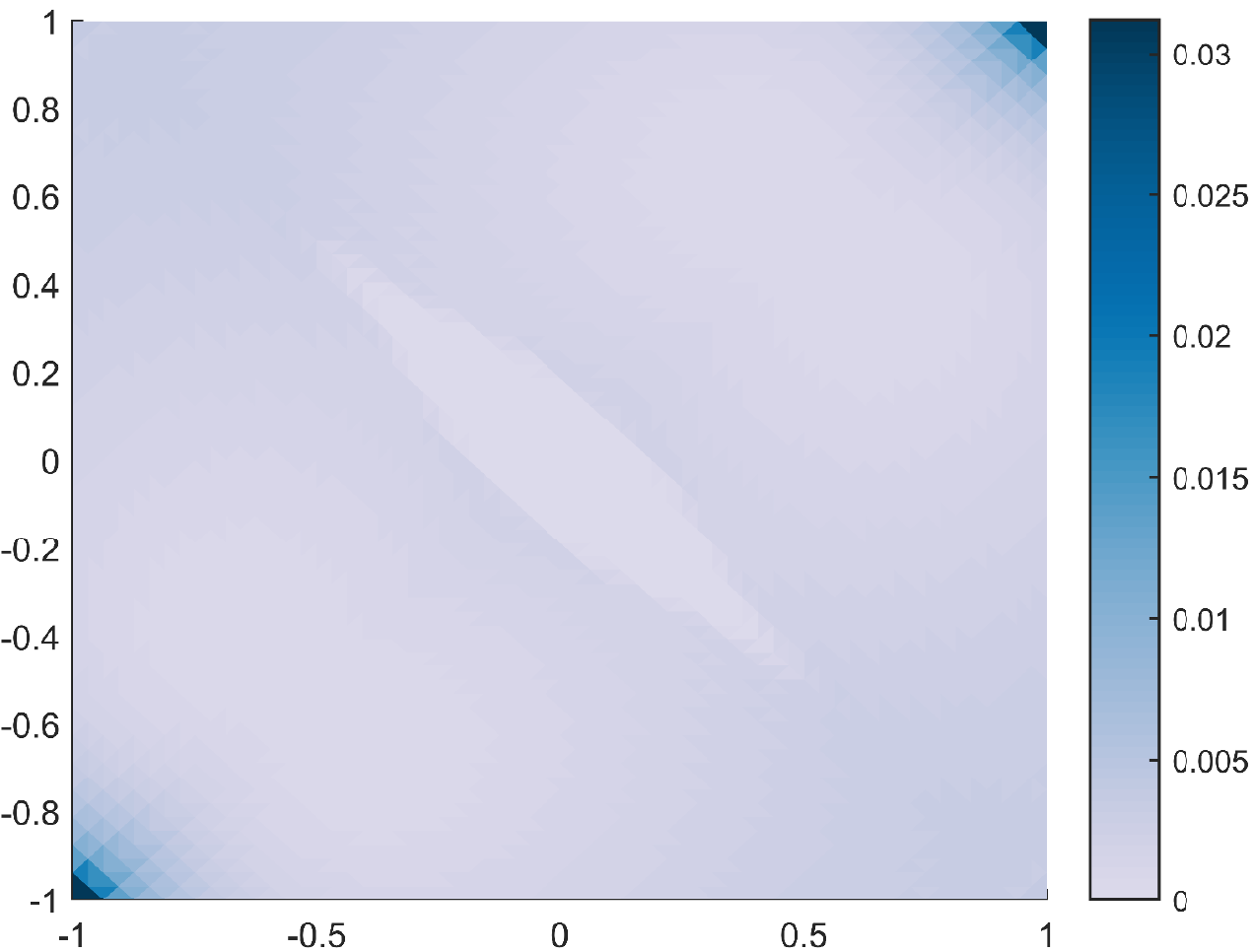}
\end{minipage}
\begin{minipage}{0.32\textwidth}
\includegraphics[width=\textwidth]{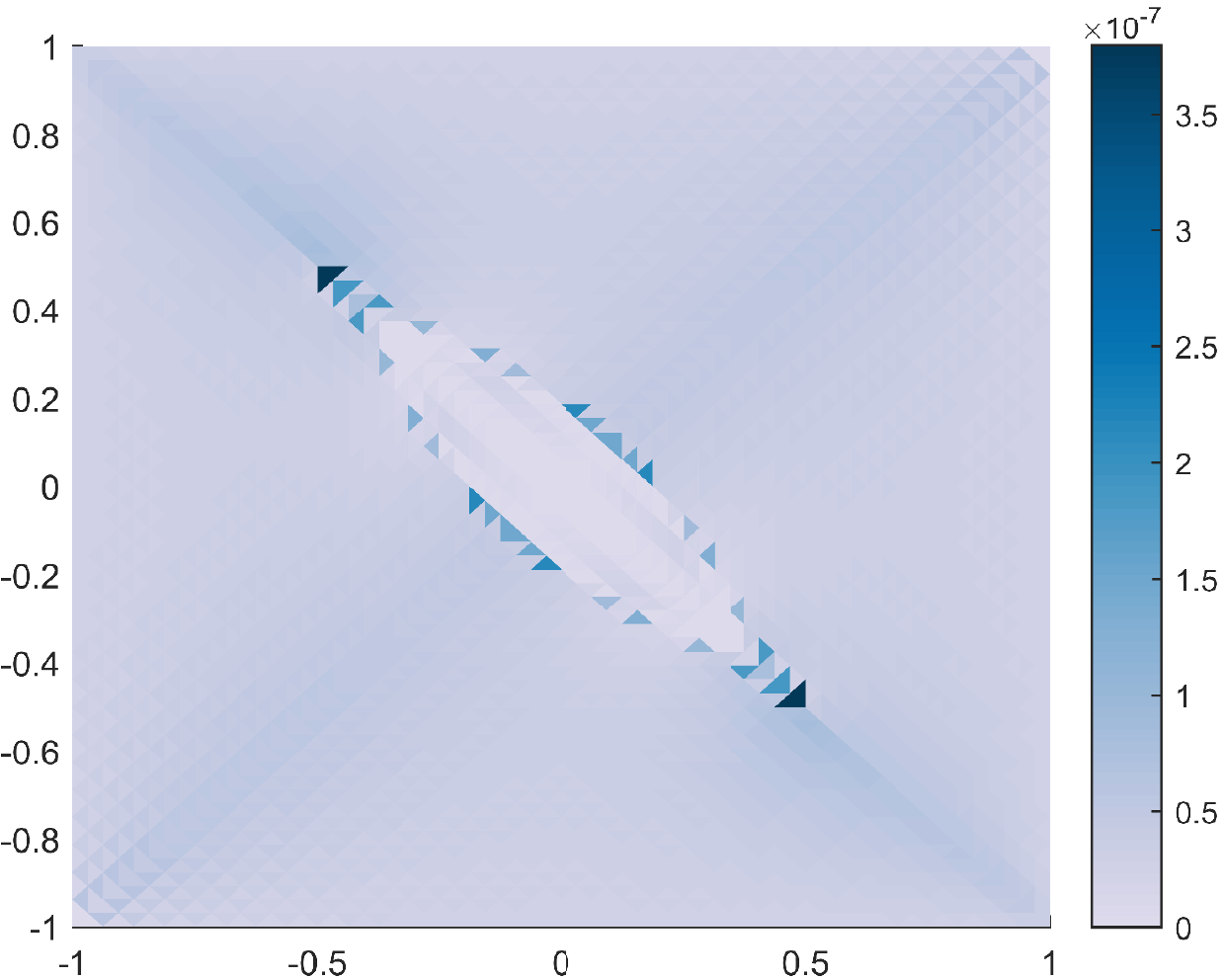}
\end{minipage}
\begin{minipage}{0.32\textwidth}
\includegraphics[width=\textwidth]{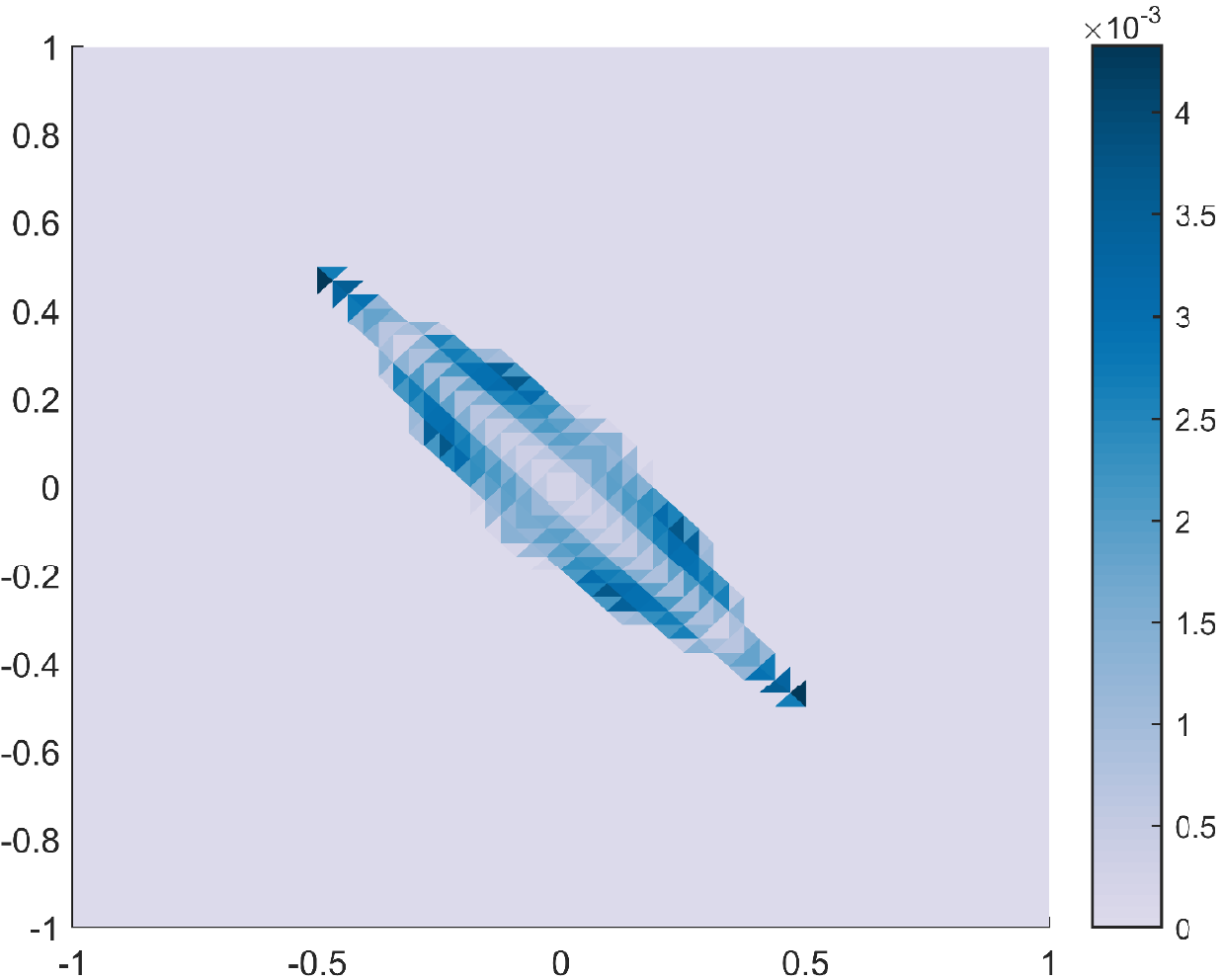}
\end{minipage}
\caption{Example II - distribution of the majorant parts  $\mathcal{M_+}_1$ (left), $\mathcal{M_+}_2$ (middle), $\mathcal{M_+}_3$ (right)
computed on level 5 triangular mesh (referred to as level 5 in Table \ref{table_exampleII}).}
\label{figure_majorantParts_exampleII}
%\vspace{1cm}
%\begin{minipage}{0.45\textwidth}
%\includegraphics[width=\textwidth]{example_II_error}
%\end{minipage}
%\hspace{1cm}
%\begin{minipage}{0.45\textwidth}
%\includegraphics[width=\textwidth]{example_II_majorant}
%\end{minipage}
%\caption{Example II - distribution of the error $\frac{1}{2}\|\nabla(u_{\lambda_h}-u_{ref})\|^{2}_{\Omega,A}$ (left) and of the majorant $\mathcal{M_+}$  (right) computed on a triangular mesh (referred to as level 5 mesh in Table \ref{table_exampleII}).}
%\label{figure_error_and_majorant_exampleII}
\end{figure}

The second example is also taken from \cite{Bo} and considers a square domain
\begin{equation}  \label{example2_domain}
\Omega=X \times Y:=(-1, 1) \times (-1, 1),
\end{equation}
constant coefficients
\begin{equation}  \label{example1_coefficients}
\alpha_{+}=\alpha_{-}=4.
\end{equation}
The Dirichlet boundary conditions as assumed in the form
\begin{equation}  \label{example2_BC}
 u(x,y)=\left\{\begin{array}{ll} x+1 \quad & x \in [-1, 1] \quad \mbox{  and  } \quad  y=1, \\
                                                   x-1 \quad & x \in [-1, 1]  \quad \mbox{  and  } \quad  y=-1, \\
 			                   y+1 \quad & y \in [-1, 1] \quad \mbox{  and  } \quad  x=1, \\
                                                   y-1 \quad & y \in [-1, 1]  \quad \mbox{  and  } \quad  x=-1.
                                    \end{array} \right.
\end{equation}
The exact solution  $u \in K$ is not known for this example. Consequently, no apriori information about the shape of the free boundary or the value of the exact energy $J(u)$ is provided. We compute approximation pairs $(\lambda_h,u_{\lambda_h}) \in \Lam_h \times K_h$ again for a sequence of nested uniformly refined meshes. Levels 1 and 2 meshes are depicted in Figure \ref{figure_meshes_exampleII}. An approximative solutions pair $(\lambda_h,u_{\lambda_h}) \in \Lambda_h \times K_h$ obtained by the dual-based solver is depicted  in Figure \ref{figure_solution_exampleII}. The approximative Lagrange multiplier field $\lambda_h$ presumably indicates the exact free boundary. Space distributions of majorant subparts are visualized in Figure \ref{figure_majorantParts_exampleII}. We assume that the density of $\mathcal{M_+}_3$ serves as an indicator of the exact free boundary. Table \ref{table_exampleII} summarizes computations on all nested uniformly refined triangular meshes. The exact energy $J(u)$ is not known but it is replaced by the energy $J(u_{ref})$ of a reference solution $u_{ref}$ in Table  \ref{table_exampleII}. The reference solution $u_{ref}$ is computed as $u_{\lambda_h}$ on the mesh one level higher (level 6 uniformly refined triangular mesh here).

%, the distributions of exact error $\frac{1}{2}\|\nabla(u_{\lambda_h}- u_{ref})\|^{2}_{\Omega,A}$ and of the majorant %$\mathcal{M_+}(u_{\lambda_h},\dots)$ in Figure \ref{figure_error_and_majorant_exampleII}.

\begin{table}[!h]
\begin{center}
\begin{small}
\begin{tabular}{|c|c|c|c|c|c|c|c|c|}
\hline
$\mbox{level}$  & $|\NN|$  & $J(u_{\lambda_h}) \, (I^*(\lambda_h))$ & $J(u_{\lambda_h}) -J(u_{ref})$
 & $\mathcal{M_+}(u_{\lambda_h},\cdot)$   & $\mathcal{M_+}_1(\cdot)$  & $\mathcal{M_+}_2(\cdot)$  & $\mathcal{M_+}_3(\cdot)$
%& $\frac{\sqrt{2 \mathcal{M}(v,\dots)}}{\|\nabla(v-u)\|_{\Omega,A}}$
\\
%   &  &   &  &  efficiency \\
\hline
1 & 13 & 13.6667 (13.6667) & 6.65e-01 & 2.41e+00 & 2.03e+00 & 3.30e-01 & 5.38e-02 \\
2 & 41 & 13.1924 (13.1924) & 1.90e-01 & 8.24e-01 & 7.90e-01 & 2.88e-02 & 5.26e-03 \\
3 & 145 & 13.0491 (13.0489) & 4.71e-02 & 2.18e-01 & 2.15e-01 & 2.20e-03 & 5.24e-04 \\
4 & 545 & 13.0137 (13.0133) & 1.17e-02 & 5.58e-02 & 5.46e-02 & 1.49e-04 & 1.03e-03 \\
5 & 2113 & 13.0045 (13.0041) & 2.50e-03 & 1.40e-02 & 1.36e-02 & 1.09e-05 & 4.36e-04 \\
\hline
6 & 8321 & 13.0020 (13.0019) &  \multicolumn{4}{ c }{\multirow{1}{*}{not evaluated} } & \\
\hline
\end{tabular}
\caption{Computations of Example II on various uniform triangular meshes. Note that
$J(u_{ref})=J(u_{\lambda_h})$ for $u_{\lambda_h}$ computed on level 6 mesh. }\label{table_exampleII}
\end{small}
\end{center}
\end{table}

\begin{figure}[th]
\center
\includegraphics[width=0.6\textwidth]{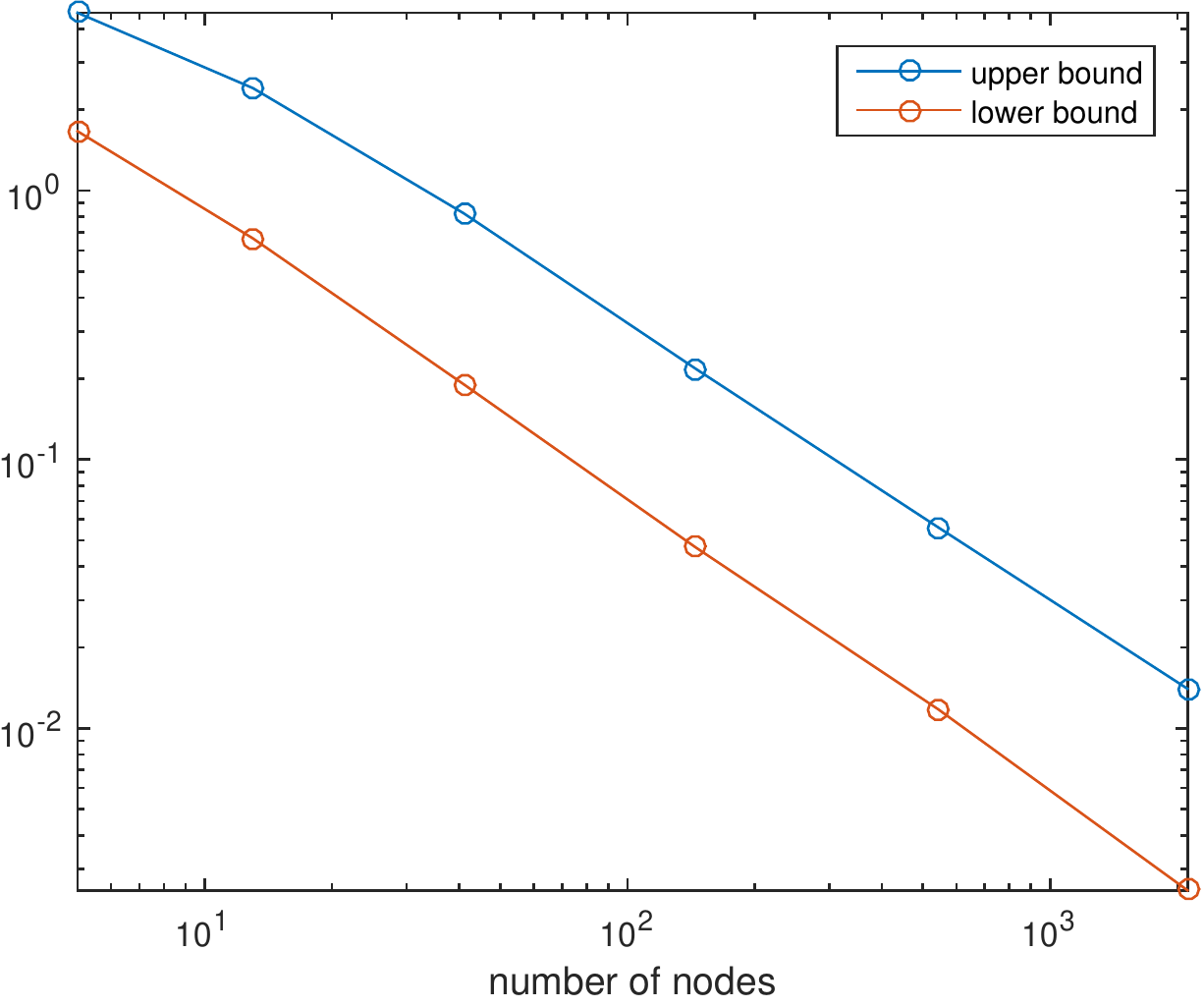}
\caption{Example II - convergence of the difference of energies $J(u_{\lambda_h}) -J(u)$ is controlled by its computable upper bound $\mathcal{M_+}(u_{\lambda_h},\dots)$ and its computable lower bound $J(u_{\lambda_h}) -J(u_{ref})$ .}
\label{figure_convergence_exampleII}
\end{figure}

\begin{remark}[Lower bound of difference of energies based on a reference solution]
If a reference solution $u_{ref}$ is available, its energy $J(u_{ref})$ satisfies $J(u_{\lambda_h}) \geq J(u_{ref}) \geq J(u)$ and
\begin{equation}  \label{both_bounds}
J(u_{\lambda_h}) -J(u_{ref}) \leq J(u_{\lambda_h}) -J(u) \leq \mathcal{M_+}(u_{\lambda_h},\dots).
\end{equation}
The inequality \eqref{both_bounds} provides actually guaranteed lower and upper bounds of the difference of energies $J(u_{\lambda_h}) -J(u).$ Figure
\ref{figure_convergence_exampleII} displays convergence of both bounds of $J(u_{\lambda_h}) -J(u)$ for considered 5 levels approximations $u_{\lambda_h}$. By refolmulating  \eqref{both_bounds} we get guranteed bounds of the exact energy
\begin{equation}\label{lower_bound_of_energy}
 J(u_{\lambda_h}) - \mathcal{M_+}(u_{\lambda_h},\dots) \leq J(u) \leq J(u_{ref}).
\end{equation}
valid for every approximation $u_{\lambda_h} \in K_h$. For this example,  lower bound of \label{lower_bound_of_energy} create an increasing sequence reported in Table \ref{table_exampleII_energy}.

\begin{table}[hb]
\begin{center}
\begin{small}
\begin{tabular}{|c|c|c|c|c|c|}
\hline
level & 1  & 2  & 3 & 4 & 5   \\
\hline
lower bound of energy $J(u)$ & 11.2582 & 12.3686 &  12.8315 &  12.9580 &  12.9905  \\
\hline
\end{tabular}
\caption{Example II: lower bound of energy $J(u)$ computed for various triangular meshes. }\label{table_exampleII_energy}
\end{small}
\end{center}
\end{table}
The sharpest available estimate of $J(u)$ based on the largest lower bound from Table    \ref{table_exampleII_energy} and  the smallest upper bound from Table \ref{table_exampleII} and  read
$$12.9905 \leq J(u) \leq 13.0020$$
and it suggests $J(u)=13$  although there is no analytical proof of it.
\end{remark}

\subsection{Implementation details} \label{subsec:software}
Both numerical examples are implemented in MATLAB and the code available  for download at
\begin{center}
\url{http://www.mathworks.com/matlabcentral/fileexchange/57232}
\end{center}
The code is based on vetorization techniques of \cite{AnVa, RaVa}. The main file 'start.m' is located in the directory 'solver\_two\_phase\_obstacle'. The following parameters can be adjusted:
\begin{itemize}
\item[] 'levels\_energy\_error' - the number of the finest uniform triangular level (default is '5')
\item[] 'iterations\_majorant' - the number of iterations of Algoritm \ref{alg:majorant_minimzation}  (default is '1000')
\end{itemize}
The dual based solver of Subsection \ref{subsec:solver} is implemented in  'optimize\_energy\_dual\_mu\_constant\_compact.m' and the underlying quadratic programming function 'quadprog' requires the optimization toolbox of MATLAB to be available. Evalulation of the primal energy $J(u_{\lambda_h})$ for a given function $u_{\lambda_h} \in K$ is done in the function
'energy'. This function is able  to provide an exact quadrature \cite{KaVa} of the energy $J(v)$ for any function $v \in K_h$, including nondifferentiable terms $\int_{\Omega} v^+ dx, \int_{\Omega} v^- dx.$

\section{Conclusions and future outlook}
A dual based solution algorithm to provide a finite element approximation of the Lagrange multiplier of the perturbed problem was described and tested on two benchmarks in 2D. The finite elements approximation of the primal minimization problem can be easily reconstructed from Lagrange multipliers by solving one linear system of equations. The quality of such approximation is measured in terms in terms of a fully computational functional majorant. A nonlinear part of the optimized functional majorant seems to work as an indicator of the free boundary. The functional majorant minimization is based on a subsequent minimization and therefore requires many iterations. We would like to speed up majorant optimization in the future.

\end{document}